\documentclass[12pt]{amsart}
\usepackage[left=1.28in,top=.8in,right=1.28in,bottom=.8in]{geometry}

\usepackage{mathptmx}
\usepackage{amsmath}
\usepackage{amssymb}
\usepackage{amsthm,amscd}
\usepackage{amssymb}
\usepackage[T1]{fontenc}
\usepackage{mathtools}
\usepackage{latexsym}
\usepackage{enumerate}
\usepackage{enumitem}
\usepackage{graphics}
\usepackage[dvipsnames]{xcolor}
\usepackage[toc,page]{appendix}
\usepackage{graphicx}
\usepackage{cite}
\usepackage[strict]{changepage}
\usepackage{hyperref}
\hypersetup{
    colorlinks=true,
    urlcolor=gray,
    citecolor=gray,   
    linkcolor=gray,
}
\usepackage{tikz-cd}
\usepackage{tikz}
\usetikzlibrary{matrix,fit,backgrounds}
\usepackage{xifthen}
\usepackage{hyperref}
\usepackage{import}
\usepackage[final]{fixme}
\usepackage[final]{changes}

\theoremstyle{plain} 
\newtheorem{theorem}{Theorem}[section]
\newtheorem{lemma}[theorem]{Lemma}
\newtheorem{proposition}[theorem]{Proposition}
\newtheorem{corollary}[theorem]{Corollary}

\newtheorem{example}[theorem]{Example}
\newtheorem{examples}[theorem]{Examples}

\newtheorem{thmx}{Theorem}

\newtheorem{propx}[thmx]{Proposition}
\newtheorem{corx}[thmx]{Corollary}
\newtheorem{theorem*}{Theorem}[]
\newtheorem{question*}[theorem*]{Question}
\newtheorem{proposition*}[theorem*]{Proposition}

\theoremstyle{definition}
\newtheorem{defn}[theorem]{Definition}
\newtheorem{remark}[theorem]{Remark}

\DeclareFontFamily{U}{wncy}{}
\DeclareFontShape{U}{wncy}{m}{n}{<->wncyr10}{}
\DeclareSymbolFont{mcy}{U}{wncy}{m}{n}
\DeclareMathSymbol{\Sha}{\mathord}{mcy}{"58} 

\DeclareMathOperator{\pr}{pr}

\DeclareMathOperator{\Mod}{Mod}

\DeclareMathOperator{\Isom}{\underline{Isom}}
\DeclareMathOperator{\id}{id}
\DeclareMathOperator{\Hom}{mor}

\DeclareMathOperator{\Aut}{Aut}

\DeclareMathOperator{\End}{End}
\DeclareMathOperator{\op}{op}

\DeclareMathOperator{\Spec}{Spec}
\DeclareMathOperator{\calO}{\mathcal{O}}
\DeclareMathOperator{\HH}{H}

\DeclareMathOperator{\calC}{\mathcal{C}}
\DeclareMathOperator{\calF}{\mathcal{F}}

\DeclareMathOperator{\Band}{Band}

\DeclareMathOperator{\SchF}{(Sch/\textit{F})_{\acute{e}t}}

\DeclareMathOperator{\calP}{\mathcal{P}}

\DeclareMathOperator{\G}{\mathcal{G}}
\DeclareMathOperator{\PP}{\mathcal{P}}
\DeclareMathOperator{\C}{\mathcal{C}}
\DeclareMathOperator{\SP}{Sp}
\DeclareMathOperator{\Tor}{TORSOR}
\DeclareMathOperator{\Bitorsor}{BITORSOR}
\DeclareMathOperator{\F}{\mathcal{F}}
\DeclareMathOperator{\V}{\mathcal{V}}
\DeclareMathOperator{\PAP}{PP}
\DeclareMathOperator{\TPAP}{TPP}
\DeclareMathOperator{\Isex}{\underline{Isex}}
\DeclareMathOperator{\ad}{ad}
\DeclareMathOperator{\obj}{obj}
\DeclareMathOperator{\BPAP}{BPP}
\DeclareMathOperator{\GPAP}{GPP}
\DeclareMathOperator{\cP}{\mathcal{P}}
\DeclareMathOperator{\invlim}{\varprojlim}

\DeclareMathOperator{\Vect}{VECT}
\DeclareMathOperator{\GL}{GL}

\DeclareMathOperator{\Gerbe}{Gerbes}
\DeclareMathOperator{\K}{\mathcal{K}}
\DeclareMathOperator{\Nrd}{Nrd}

\DeclareMathOperator{\YY}{\mathcal{Y}}
\DeclareMathOperator{\XX}{\mathcal{X}}
\DeclareMathOperator{\calB}{\mathcal{B}}
\DeclareMathOperator{\calU}{\mathcal{U}}
\DeclareMathOperator{\Xhat}{\tilde{X}}

\DeclareMathOperator{\SL}{SL}
\DeclareMathOperator{\PGL}{PGL}
\DeclareMathOperator{\Z}{\mathbb{Z}}

\DeclareMathOperator{\chara}{char}
\DeclareMathOperator{\sep}{sep}

\newcommand{\Sch}[1]{(\text{Sch} / #1)}

\fxusetheme{color}
\definechangesauthor[color=orange]{Bastian}
\title{Gerbe patching and a Mayer-Vietoris sequence over arithmetic curves} 
\author{Bastian Haase}
\address{Department of Mathematics \& Computer Science, Emory University, 400 Dowman Drive
NE, Atlanta, GA 30322, USA}
\email{bhaase@emory.edu}

\begin{document}

\maketitle

{\hypersetup{
   linkcolor=black,
}

\listoffixmes

\begin{abstract}
  We discuss patching techniques and local-global principles for gerbes over arithmetic curves. Our patching setup is that
  introduced by Harbater, Hartmann and Krashen (cf. \cite{Harbater2010},\cite{Harbater2011}).
  Our results for gerbes
  can be viewed as a 2-categorical analogue on their results for torsors. Along the way, we also discuss
  bitorsor patching and local-global principles for bitorsors.
  As an application of these results, we obtain a Mayer-Vietoris sequence with respect to patches for non-abelian 
  hypercohomology sets with values in the crossed module $G \rightarrow \Aut(G)$ for
  $G$ a linear algebraic group. 
  Using local-global principles for gerbes, we also prove local-global principles for points on 
  homogeneous spaces under linear algebraic groups $H$ that are special (e.g. $\SL_n$ and $\SP_{2n}$)
  for certain kind of stabilizers.
\end{abstract}

\section{Introduction}
\let\thefootnote\relax\footnote{\textbf{2010 Mathematics Subject Classification:} Primary: 14H25, 18G50. Secondary: 14M17. \newline 
\textbf{Key words and phrases:} gerbes, non-abelian cohomology, local-global,
patching, Mayer-Vietoris, homogeneous spaces.}
Starting with the famous Hasse principle for quadratic forms over number fields,
local-global principles for algebro-geometric objects over number fields
have been studied extensively. In \cite{Harbater2010}, Harbater, Hartmann and Krashen introduced
patching techniques for function fields over arithmetic curves. These patching techniques in turn
allowed the authors to deduce local-global principles for quadratic forms and central simple algebras
over these  fields.

Given a function field of an arithmetic curve, the patching setup consists of a finite set of fields $F_i$ over
$F$. The inverse limit of these fields with respect to certain inclusions is $F$.
The definition of the fields $F_i$ is deduced from a geometric ''covering'' of the special fiber
of the curve. Let $\calF$ denote the system of fields $\{F_i \}$. 
When we say that patching holds for a certain class of objects over $\F$, we mean that the datum
of such an object over $F$ is equivalent to the datum of a collection of compatible
objects over the fields $F_i$. We say that the local global principle holds if triviality
of an object over $F$ is equivalent to triviality on the level of $F_i$ for all $i$.

In \cite{Harbater2012a} and  \cite{Harbater2011}, Harbater, Hartmann and Krashen
prove patching results for torsors under linear algebraic groups and Galois cohomology groups.
 In both of these papers,
the authors use their patching results to characterize when local-global principles hold.

In this paper, we extend these results to the case of gerbes over arithmetic curves.
We prove patching results for certain gerbes and  use these results to give
a criterion for when a local-global principle for gerbes hold. Using this criterion,
we then prove that the local-global principle holds for gerbes banded by certain linear algebraic groups.

Note that for a sheaf of abelian groups $A$ defined over $F$, its second cohomology group
classifies gerbes banded by $A$. 

Gerbes can not only be thought of as a 2-categorical analogue of torsors, but also arise naturally
in the context of algebraic geometry. For instance, if $X$ is a homogeneous space under a group sheaf $H$,
then the quotient stack $[X/H]$ that describes relative morphisms of principal homogeneous spaces under $H$
to $X$ is a gerbe (whose band is determined by the stabilizers of geometric points of $X$).

Our main result on gerbe patching is the following theorem.
\begin{thmx}[cf. Theorem \ref{curve-gerbes}]
  Let $L$ be a linear algebraic band over $F$.
  Assume either  
  \begin{itemize}
    \item $\chara(k)=0$  or
    \item $\chara(k)=p>0$ and  the center of $L$ has finite order not divisible by $\chara(k)$ and
      $G$-bitorsor factorization holds for $\F$ over   every cover $Z \rightarrow Y$ .
  \end{itemize}
  Then, gerbe patching holds over $\F$ for $L$-banded gerbes.
\end{thmx} 
Our approach to gerbe patching relies on a semi-cocyclic description of gerbes in terms of
bitorsors introduced by Breen in \cite{Breen1990} and \cite{Breen1991}.

Recall that a $(G,H)$-bitorsor $P$ is a left $G$- and right $H$-torsors such
that the left $G-$ and the right $H$-action commute. In the case $G=H$ we simply
speak of a $G$-bitorsor. The trivial $G$-bitorsor is $G$ with the natural right and left action.

In view of the semi-cocyclic description, we  prove gerbe patching by proving
that patching holds for bitorsors.
 Our results on patching gerbes and bitorsors allow a natural
interpretation in terms of hypercohomology with values in the crossed module
$G \rightarrow \Aut(G)$.
\begin{corx}[Mayer-Vietoris of non-abelian hypercohomology over curves, Corollary \ref{mayer-curves}]
  Let $G$ be a linear algebraic group defined over $F$ and let $L$ denote the
  associated band.
  Under the assumption of Theorem A, there is an exact sequence
  of pointed sets
   \begin{center}
    \begin{tikzcd}[row sep =small, column sep = small]
       \HH^{-1}(F, G \rightarrow \Aut(G)) \arrow{r}&
    \prod_{i \in I_v} \HH^{-1}(F_i, G \rightarrow \Aut(G))\arrow{r} \arrow[draw=none]{d}[name=Z, shape=coordinate]{} & 
    \prod_{k \in I_e} \HH^{-1}(F_k, G \rightarrow \Aut(G))  \arrow[rounded corners,to path={ -- ([xshift=2ex]\tikztostart.east)|- (Z) [near end]\tikztonodes-| ([xshift=-2ex]\tikztotarget.west)-- (\tikztotarget)}]{dll} \\
      \HH^0(F, G \rightarrow \Aut(G)) \arrow{r}&
    \prod_{i \in I_v} \HH^0(F_i, G \rightarrow \Aut(G))\arrow{r} \arrow[draw=none]{d}[name=Z, shape=coordinate]{} & 
    \prod_{k \in I_e} \HH^0(F_k, G \rightarrow \Aut(G))  \arrow[rounded corners,to path={ -- ([xshift=2ex]\tikztostart.east)|- (Z) [near end]\tikztonodes-| ([xshift=-2ex]\tikztotarget.west)-- (\tikztotarget)}]{dll} \\
      \HH^1(F, G \rightarrow \Aut(G)) \arrow{r}& 
    \prod_{i \in I_v} \HH^1(F_i, G \rightarrow \Aut(G))\arrow[yshift=3]{r} \arrow[yshift=-3]{r}&
    \prod_{k \in I_e} \HH^1(F_k, G \rightarrow \Aut(G))
    \end{tikzcd}
  \end{center}
\end{corx} 
Note that this result is analogous to the Mayer-Vietoris sequences obtained in \cite{Harbater2012a} and \cite{Harbater2011}.

By use of the Mayer-Vietoris sequence, we can prove that the local-global principle for $G$-bitorsors holds if and only if
the center $Z(G)$ of $G$ satisfies simultaneous factorization over $\F$ (see Section \ref{patching} for details). In particular, using
a result of Harbater, Hartmann and Krashen in \cite{Harbater2009}, we prove that local-global principle for
$G$-bitorsors holds if $Z(G)$ is rational and connected.

The Mayer-Vietoris sequence also allows us to prove that the local-global principle for $G$-gerbes is equivalent
to factorization for $G$-bitorsors (see Section \ref{generalities-bitorsors} for details).
We prove that bitorsor factorization (and thus local-global principle) 
holds for various groups such as $\SL_1(D)$ where $D$ is a central simple algebra over $F$.

As an application of the theory, we use local-global principles for gerbes to deduce local-global
principles for certain homogeneous spaces. Recall that a linear algebraic group $H$ over $F$ is
\emph{special} if  $\HH^1(K,H)$ is trivial for every field extension $K/F$. The groups
$\SL_n$ and $\SP_{2n}$ are examples of special groups. Let now $X$ be a homogeneous space
under a special group $H$ and assume that the stabilizer $G$ is defined over $F$. 
\begin{thmx}[cf. Theorem \ref{homspace}]
  Let $H$ be special and let $X$ be homogeneous space under $H$. Assume that the stabilizer $G$ of a geometric point of $X$ is defined
  over $F$. 
  If $\chara(k)=0$, the local-global principle for $X$ holds if and only if
  bitorsor factorization holds for $G$. \par
  If $\chara(k)=p>0$,  assume that 
  \begin{itemize}
    \item $p \nmid |Z(G)| < \infty$,
    \item $G$-bitorsor factorization holds for $\F$ over every cover $Z \rightarrow F$.
  \end{itemize} 
  Then, the local-global principle holds for $X$ over $\F$.
\end{thmx} 
Using our results on bitorsor factorization, we then show that local-global principle
for $X$ holds for certain type of $G$, e.g. if $G$ is finite and simple.

The paper is structured as follows.
In the second section we review the patching setup introduced by Harbater, Hartmann and Krashen and describe
the main results concerning patching over arithmetic curves.

In the following two sections, we discuss bitorsor and gerbe patching over an arbitrary inverse factorization
system of fields and relate it to patching of torsors.

In the third section, we utilize a description of bitorsors introduced by Breen in \cite{Breen1990}
to reduce bitorsor patching to torsor patching. This then allows us to deduce the first two rows of 
the Mayer-Vietoris sequence in hypercohomology. From this sequence, we can then show that
local-global principle for $G$-bitorsors with respect to patches holds if and only if the
center of $G$ satisfies factorization.

In the fourth section, we recall a semi-cocyclic description of gerbes in terms of bitorsors.
We also recall the notion of a band and prove the following patching results for its second non-abelian
cohomology set.
\begin{propx} \ref{assum-gerbe-charp}
  Let $L$ be a band over $F$ such that its center $Z(L)$ is a linear algebraic group over $F$
  with finite order coprime to $\chara(k)$. Then, if $\HH^2(F,L) \neq \emptyset$, 
  patching holds for $\HH^2(\bullet, L)$ over
  $\F$.
\end{propx} 
Using our results on bitorsor patching, we can  prove that gerbe patching holds
 under a further technical assumption (cf. Theorem \ref{gerbe-patching}).
This allows us in turn to complete the construction of the Mayer-Vietoris sequence above.
Finally, we relate the local-global principle for gerbes to factorization of bitorsors.

In the fifth section, we then specialize to patching over arithmetic curves  and prove our main results for patching bitorsors
and gerbes in this setting. We also prove that bitorsor factorization and 
local-global principle for gerbes hold for various interesting linear algebraic groups.
We conclude the section by applying our results on local-global principles for gerbes
to deduce our results on local-global principles for points on  homogeneous spaces.

In the appendix, we review hypercohomology with values in crossed modules as developed by
Breen in \cite{Breen1990} and Borovoi in \cite{Borovoi1992}.

\subsection{Acknowledgments}
The author would like to thank  Max Lieblich for various helpful conversations.
The author is partially supported by National
Science Foundation grants DMS-1401319 and DMS-1463882.
\section{Patching and Local-Global Principle}\label{patching}

\subsection{Patching of vector spaces}
In this section, we will recall patching results from Harbater, Hartman and Krashen
(compare \cite{Harbater2010},\cite{Harbater2011} and \cite{Harbater2014}).
We will adapt their notation.

Throughout this section, let $\F= \left \{F_i \right \}_{i \in I}$  denote a finite inverse systems of fields
with inclusions as morphisms. Let $F$ denote $\invlim \F$. All sheaves will be sheaves in the big \'{e}tale topology
over $\Sch{F}$ or $\Sch{F_i}$ for $i \in I$.

\begin{defn}
  A \emph{factorization inverse system} over a field $F$ is a finite inverse system of fields such that 
  \begin{enumerate}
    \item $F$ is the inverse limit.
    \item The index set $I$ can be partitioned as $I=I_v \sqcup I_e$ such that: 
      \begin{enumerate}
        \item For any $k \in I_e$ there are exactly two elements $i,j \in I_v$ such that $i,j >k$.
        \item These are the only relations.
      \end{enumerate}
  \end{enumerate}
\end{defn}

For each index $ k \in I_e$, fix a labeling $l_k,r_k$ for the two elements
in $I_e$ with $l_k,r_k>k$. Then, let $S_I$ denote the set of triples
$(l_k,r_k,k)$.

Given a factorization inverse system, one can associate to it a (multi-)graph $\Gamma$. Its vertices are the elements of $I_v$, where the edges come from $I_e$ (explaining the subscripts). 
Given $k\in I_e$, the corresponding edge connects the vertices $i,j \in I_v$ iff $(i,j,k) \in S_I$.
Note that  $\Gamma$ is connected, as the inverse limit F would otherwise admit zero-divisors. 
We will sometimes specialize to the case where $\Gamma$ is a tree.

\begin{example}
  A basic example of a factorization inverse system is given by fields 
  $F\subset F_1, F_2 \subset F_0 $ such that $F=F_1 \cap F_2$. Pictorially, we get
  
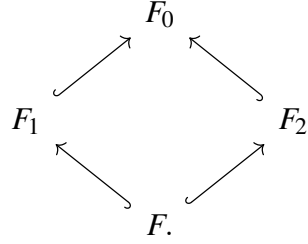
\begin{center}
 \begin{tikzcd}
   \; & F_0 & \; \\
   F_1 \arrow[hook]{ru} & \; & F_2 \arrow[hook]{lu} \\
   \; & F.  \arrow[hook]{ru} \arrow[hook]{lu}& \; 
\end{tikzcd}
\end{center}
Note that in this case $I_e=\{0 \}, I_v=\{1,2\}$ and $S_I=\{(1,2,0)\}$. The corresponding
graph $\Gamma$ is a tree.
\end{example}

\begin{defn}
  A vector space patching problem $\V=\left(\left \{V_i \right\}_{i \in I_v}, \{\nu_k\}_{k \in I_e}\right)$ for
  a factorization inverse system $\F$ is given by finite dimensional $F_i$ vector spaces
  $V_i$ together with $F_k$ isomorphisms $ \nu_{k}\colon V_i \otimes_{F_i}F_k \rightarrow V_j \otimes_{F_j} F_k$
  whenever $(i,j,k) \in S_I$.  \par
  A morphism of patching problems $\V \rightarrow \V'$ is a collection of morphisms
  $ V_i \rightarrow V_i'$ for all $i \in I_v$ that are compatible with the morphisms
  $\nu_{k}, \nu_{k}'$. \par
  The category of vector space patching problems is denoted by $\PAP(\F)$. By construction, 
  we have $\PAP(\F) \simeq \prod_{(i,j,k) \in S_I} \Vect(F_i) \times_{\Vect(F_k)} \Vect(F_j)$.
  If $A/F$ is a finite product of finite separable field extension, let $\PAP(\F_A)$ denote the category of free module patching problems: objects are collections 
  $\left(\{M_i\}_{i \in I_v},\{\nu_k \}_{k \in I_e}\right)$ where $M_i$ is a free module over $A_i=F_i\otimes_F A$ of finite rank and $\nu_k \colon M_i|_{A_k} \rightarrow M_j|_{A_k}$ are isomorphisms of $A_k$ modules.
\end{defn}

Note that we have a canonical functor 
\begin{align*}
  \beta \colon \Vect(F) \rightarrow \PAP(\F)
\end{align*}
where $\Vect(F)$ is the category of finite dimensional vector spaces over $F$.

\begin{defn}
  A \emph{solution} to a vector space patching problem $\V$ is a $F$ vector space $V$
  such that $\beta(V)$ is isomorphic to $\V$. \par
\end{defn}

If $\beta$ is an equivalence of categories, then every patching problem has a solution that
is unique up to isomorphism.

\begin{defn}
  A \emph{linear algebraic group} $G$ over a field $K$ is a smooth affine group scheme of finite type.
\end{defn}
\begin{defn}
  A linear algebraic group $G$ satisfies \emph{simultaneous factorization over $\F$}
  if for any collection of elements $a_k \in G(F_k)$ with $k \in I_e$, there are elements
  $a_i \in G(F_i)$ for all $i \in I_v$ such that $a_k=a^{-1}_r a_l \in G(F_k)$ for all
  $(l,r,k) \in S_I$. \par
  In the case where $G=\GL_n$, we simply say that \emph{simultaneous factorization holds over $\F$}.
\end{defn}

The following result provides the basis of many patching results obtained by Harbater, Hartman
and Krashen. 
\begin{proposition}[{\cite[Proposition 2.1]{Harbater2010}}] \label{matrix}
The functor $\beta \colon \Vect(F) \rightarrow \PAP(\F)$ is an equivalence
of categories if and only if simultaneous factorization holds over $\F$.  
\end{proposition}

If $A=\prod_{i=1}^{n} L_i$ is a finite product of finite separable field extensions of $F$, then we denote by
$\F_A$ the inverse system $\{A_i:=A\otimes_F F_i \}_{i \in I}$. It is not necessarily an inverse factorization
system 
but we have $A = \varprojlim \F_{A}$.
Let $\Mod(A)$ denote the category of free modules of finite rank over $A$. Let $\PAP(\F_A)$ denote
the category of patching problems of the form $\left(\{M_i \}_{i \in I_v}, \{\nu_k \}_{k \in I_e}\right)$ where 
$M_i$ is a free $A_i$-module of finite rank and $\nu_k \colon M_i|_{A_k} \rightarrow M_j|_{A_k}$
is an isomorphism of $A_k$-modules for $(i,j,k) \in S_I$. Note that we have a natural functor
$\widehat{\beta} \colon \Mod(A) \rightarrow \PP(\F_A)$.

The following proposition is a slight variant of \cite[Lemma 2.2.7]{Harbater2014}.
\begin{proposition} \label{vect}
Assume that patching holds over $\F$, i.e. that the functor $\beta \colon \Vect(F) \rightarrow \PAP(\F)$ is an equivalence.
Let $A=\prod_{i=1}^nL_i$ be a product of finitely many finite separable field extensions.
 Then, patching holds over $\F_A$, i.e.  the natural  functor
$\widehat{\beta} \colon \Mod(A) \rightarrow \PAP(\F_A)$ is an equivalence.
\end{proposition} 
\begin{proof}
  Given a patching problem  $\left(\{M_i\}_{i \in I_v},\{\nu_k \}_{k \in I_e}\right)$ in $\PAP(\F_A)$, note that
  $M_i$ is a finite dimensional vector space over $F_i$ for all $i \in I_v$. Also, $\nu_k$ is an isomorphism
  of $F_k$ vector spaces. Hence, by assumption, there is an $F$-vector space $M$ together with $F$-vector space isomorphisms
  $\phi_i \colon M|_{F_i} \rightarrow M_i$ for all $i \in I_v$ that are compatible with $\nu_k$. The $A_i$-module structure 
  of $M_i$ is equivalent to the datum of  a morphism $\alpha_i \colon A_i \rightarrow \End_{F_i}(M_i)$ of $F_i$-vector spaces. As the $\nu_k$ are $A_k$ module
  morphisms, they are compatible with $\alpha_i$. Hence, as the functor $\beta$ is full, there is a morphism $\alpha \colon A \rightarrow \End_F(M)$
  of $F$ vector spaces. The resulting $A$-module $M$ solves the patching problem. This shows that $\widehat{\beta}$ is essentially surjective.
  As every morphism of $A_i$-modules is in particular a morphism of $F_i$-vector spaces, it follows that $\widehat{\beta}$ is faithful.
  Given two $A$-modules $M$, $N$ and a morphism $\widehat{\beta}(M) \rightarrow \widehat{\beta}(N)$, note that we can lift it to an $F$-vector space
  morphism $\gamma \colon M \rightarrow N$. As the the image of this morphism in $\PAP(\F_A)$ commutes with the $A_i$-action and as $\widehat{\beta}$ is faithful,
  it follows that $\gamma$ is an $A$-module morphism.
\end{proof}

\subsection{Patching of torsors and Galois cohomology} \label{gal}
Fix a field $F$ and an inverse factorization system $\F$ over $F$. Let $G$ be a group sheaf in the big \'{e}tale site over $F$.

\begin{defn}
  A \emph{left principal homogeneous space}  $T$ over $F$ under a $F$-group scheme $G$ is a $F$-scheme $T$ together with a left $G$ action
  $G \times_{F} T \rightarrow T$ such that the induced morphism 
  \begin{align*}
    G \times_{F} T \rightarrow T \times_F T
  \end{align*}
  given by $(g,t) \mapsto (g \cdot t,t)$
  is an isomorphism. \par
  A \emph{morphism of left principal $G$-homogeneous spaces} $T,T'$  over $F$ is a $F$-morphism $T \rightarrow T'$ that is $G$-equivariant.
\end{defn}

\begin{remark}
  When we speak of a principal homogeneous space, we mean \emph{left} principal homogeneous space unless explicitly stated otherwise.
\end{remark}

\begin{defn}
  A \emph{left \'{e}tale $G$-torsor} $T$ over $F$ is an \'{e}tale sheaf  with a left action of the sheaf $h_G$  such that: 
  \begin{enumerate}
    \item For all $K$-schemes $Y$ there is an \'etale cover $\{Y_i \rightarrow Y \}$ such that $T(Y_i) \neq \emptyset$.
    \item The map $h_G \times T \rightarrow T \times T$ given by $(g,t) \mapsto (g \cdot t, t)$ is an isomorphism.
  \end{enumerate}
  A \emph{morphism of  $G$-torsors} $T,T'$ over $K$ is a morphism of sheaves $T \rightarrow T'$ which is $G$-equivariant.
\end{defn}

\begin{remark}
  If $G$ is affine over the base field, then the notions of (left) principal homogeneous space and (left) torsor coincide, i.e. the category
  of (left) principal $G$-homogeneous spaces over $F$ is equivalent to the category of (left) $G$ torsors over $F$ (compare Proposition 4.5.6 in \cite{Olsson2016}).
\end{remark}

Let $\Tor(G)(F)$ denote the category of $G$-torsors over $F$. Let $\TPAP(G)(\F)$ denote the category of torsor patching problems,
i.e. the category with objects $\left(\left \{T_i \right\}_{i \in I_v}, \{\nu_k\}_{k \in I_e}\right)$ where $T_i$ is 
a $G_i$-torsors over $F_i$ and $\nu_k \colon T_i|_{F_k} \rightarrow T_j|_{F_k}$ is an isomorphism of $G_k$-torsors
for $(i,j,k) \in S_I$. A morphism of torsor patching problems
\begin{align*}
  \left(\left \{T_i \right\}_{i \in I_v}, \{\nu_k\}_{k \in I_e}\right) \rightarrow \left(\left \{T'_i \right\}_{i \in I_v}, \{\nu'_k\}_{k \in I_e}\right)
\end{align*}
is  a collection of $G_i$-torsor morphisms $T_i \rightarrow T_i'$ compatible with $\nu_k$ and $\nu_k'$.

As in the case of vector spaces, there is a natural functor 
\begin{align*}
  \beta_G' \colon \Tor(G)(F) \rightarrow \TPAP(G)(\F).
\end{align*}
where $\TPAP(G)(\F)$ is defined analogously to $\PAP(\F)$.

\begin{theorem}[{\cite[Theorem 2.3]{Harbater2011}}] \label{tors}
Let $G$ be a linear algebraic group over $F$. If the natural functor
$\beta \colon \Vect(F) \rightarrow \PAP(\F)$ is an equivalence of categories,
then so is the functor $\beta_G' \colon \Tor(G)(F) \rightarrow \TPAP(G)(\F)$.  
\end{theorem}

 The following proposition is a slight variant of \cite[Theorem 2.2.4(c)(iii)]{Harbater2014}.
\begin{proposition} \label{trs}
  Let $A=\prod_{r=1}^{n} L_r$ be a product of finitely many finite separable field extensions $L_i/F$.
  Let $G$ be a linear algebraic group over $A$ (i.e. $G|_{L_{r}}$ is a linear algebraic group over $L_i$).
  If the functor $\beta \colon \Vect(F) \rightarrow \PAP(\F)$ is an equivalence of categories,
  then so is the functor  $\widehat{\beta}_G' \colon \Tor(G)(A) \rightarrow \TPAP(G)(\F_{A})$. 
\end{proposition}
\begin{proof}
  By  Proposition \ref{vect}, patching holds for free modules of finite rank over $\F_A$.
  Hence, $\GL_n$ satisfies factorization over $\F_A$, cf. Theorem \ref{matrix}. The proof is now
  verbatim to the proof of \cite[Theorem 2.3]{Harbater2011}.
\end{proof}

Let  $G$ be an abelian algebraic group over $F$ and let $\HH^n(F,G)$ denote the $n$-th Galois cohomology group.
For every $(i,j,k) \in S_I$, we have a map $\HH^2(F_i,G) \times \HH^2(F_j,G) \rightarrow \HH^2(F_k,G)$ given by
$(\alpha_i,\alpha_j) \mapsto \alpha_i|_{F_k} \alpha_j|_{F_k}^{-1}$. The collection of these maps gives a map
$\prod_{i \in I_v} \HH^2(F_i,G)  \rightarrow \prod_{k \in I_e} \HH^2(F_k,G)$.

\begin{defn}
  We say that \emph{patching holds for $\HH^2(\circ,G)$} over $\F$ if the sequence
  \begin{align*}
    \HH^2(F,G) \rightarrow \prod_{i \in I_v} \HH^2(F_i,G)  \rightarrow \prod_{k \in I_e} \HH^2(F_k,G)
  \end{align*}
  is exact.
\end{defn}

\subsection{Patching over arithmetic curves} \label{arithmetic}

We will now recall a patching setup for function fields of arithmetic curves obtained by Harbater, Hartmann and
Krashen in \cite{Harbater2012a}.
Following the notation of \cite{Harbater2012a}, let $T$ be a complete discretely valued ring with
field of fraction $K$, uniformizer $t$ and residue field $k$.
Let $\Xhat$ be a projective, integral and normal $T$-curve with function field $F$ and let $X$ denote its closed fiber.

For any closed point $p \in X$, let $\widehat{\calO}_{\Xhat,p}$ denote the completion of the local
ring $\calO_{\Xhat,p}$ at its maximal ideal and let $F_p$ denote the fraction field of $\widehat{\calO}_{\Xhat,p}$.
For a subset $U \subset X$, that is contained in an irreducible component of X and does not meet other components,
let $R_U$ denote the subring of $F$ of rational functions regular on $U$. Let $\widehat{R}_U$  denote its $t$-adic completion
 and $F_U$ denote the fraction field of $\widehat{R}_U$.
 For each branch of X at a closed point P, i.e., for each height one prime $b$ of $\widehat{\calO}_{\Xhat,p}$   that contains $t$,
let $\widehat{R}_b$ be the completion of $\widehat{\calO}_{\Xhat,p}$ at $b$ and let $F_b$ denote its fraction field.
 
Let $\calP\subset X$ be a non-empty set of closed points of $X$ including all points where distinct irreducible components
 of $X$ meet and all closed points where $X$ is not unibranched.
This implies that $X \setminus \calP$ is a disjoint union of finitely many irreducible affine $k$ curves. The set of these curves
will be denoted by $\calU$.

Let $p \in \calP$ and $U \in \calU$ be chosen such that $p$ is contained in the closure of $U$. 
Then, the ideal defining $U$ induces an ideal in $\calO_{\Xhat,p}$.
The branches of $U$ at $p$ are the height one prime ideals in $R_p$ containing said induced ideal.
 Let $\calB$ denote the set of all branches.
 Note that we have inclusions $F \subset F_u, F_p \subset F_b$,
 whenever $b$ is a branch
corresponding to $U$ and $p$.

These fields now form a finite inverse factorization system $\F$. With the notation from Section \ref{patching}, we have
$I_v = \calP\sqcup \calU$ and $I_e =\calB$. 

The system just described allows patching for $G$-torsors if $G$ is a linear algebraic group over $F$.
\begin{theorem}[{\cite[Theorem 6.4]{Harbater2010},  \cite[Theorem 2.3]{Harbater2011}, \cite[Proposition 3.2.1]{Harbater2014}, \cite[Theorem 3.1.1]{Harbater2012a}) \label{curve-torsor}}]
  The inverse limit of $\F$ is $F$. Furthermore, patching holds for finite dimensional
  vector spaces over $\F$, every linear algebraic group $G$ over $F$ is separably
  factorizable over  $\F$ and  $G$-torsor patching holds over $\F$.
\end{theorem}

Once patching for torsors hold, one can describe local-global principle for torsors in terms
of simultaneous factorization. We say that $G$-torsors satisfy the local-global principle
with respect to $\F$ if for any $G$-torsor $P$, we have that $P_i\simeq G_i$ for all $i \in I_v$ if
 and only if $P\simeq G$ (here, $G$ denotes the trivial $G$ torsor).

\begin{theorem}[{\cite[Theorem 3.5]{Harbater2011}}]
  Let $G$ be a linear algebraic group. Then, local-global principle for $G$-torsors holds
  if and only if $G$ satisfies simultaneous factorization over $\F$.
\end{theorem}

Let us also  recall a result describing when rational linear algebraic groups satisfy factorization.
\begin{theorem}[{\cite[Corollary 6.5]{Harbater2011}}] \label{simfact}
  Let $G$ be a rational linear algebraic group. Then, $G$ satisfies simultaneous factorization
  over $\F$ if and only if $G$ is connected or $\Gamma$ is a tree.
\end{theorem}

\begin{theorem}[{\cite[Theorem 3.1.3]{Harbater2012a}}]
  Let $G$ be an abelian linear algebraic group. If $\chara(k)=p>0$, assume furthermore that
  $p \nmid |G| <\infty$. Then, for any $n\geq 0$, patching holds for $\HH^n(\circ,G)$ over $\F$.
\end{theorem}

\section{Bitorsors}\label{bitorsor}

We will later see that bitorsor patching naturally occurs in the context of gerbe patching. Thus, in this section,
we will prove that bitorsor patching holds whenever torsor patching does.

\subsection{Generalities} \label{generalities-bitorsors}

Before we define bitorsors, let us fix some notation.
If $Y$ and $X$ are $F$-schemes for some field $F$, let $X_Y=X \times_{\Spec(F)}Y$. We will also often
write $F$ instead of $\Spec(F)$.

Even though most of what we discuss here holds over any site, we will restrict ourselves to the special case
of the big \'{e}tale site of $F$. So, all sheaves in this section are sheaves over $F$ in the \'{e}tale topology.

\begin{defn}
  Let $G,H$ be group sheaves. Then, a \emph{$(G,H)$-bitorsor} $T$ over $F$ is a sheaf over $F$ with a left action of $G$ making it
 a left $G$-torsor and a right $H$ action making it a right $H$-torsor such that these actions commute. \par
 Let $\Bitorsor(G,H)(K)$ denote the category of $(G,H)$-bitorsors over $F$.
\end{defn}

\begin{defn}
  Let $P$ be a $(G,H)$-bitorsor and let $P'$ be a $(H,G')$-bitorsor for group sheaves
  $G,G',H$. Then, we define a $(G,G')$-bitorsor $P \wedge^H P'$ as the sheafification
  of the presheaf $U \mapsto P(U) \times P'(U) / \sim$, where $(p,p') \sim (ph,hp')$
  for $p \in P(U), p' \in P'(U)$ and $h \in H(U)$. It inherits its $G$ action from
  $P$ and its $G'$ action from $P'$. We say that $P \wedge^H P'$ is the \emph{wedge product}
  of $P$ and $P'$.
\end{defn}

It follows that the wedge product defines a group structure on the set of isomorphism
classes of $G$-bitorsors.
The identity element is given by the class of $G$ and the inverse of the class of a  $G$-bitorsor $P$
is given by the class of $P^{\op}$.
We will later relate local-global principles of gerbes to a factorization 
 of bitorsors over a finite inverse
factorization system $\F$.
\begin{defn}
  Let $\F=\{F_i \}_{i \in I}$ be a finite inverse factorization system and let $G$ be a group sheaf
  over $F= \varinjlim \F$. We say that \emph{$G$-bitorsor factorization} holds over $\F$ if
  for any collection $\{P_K \}_{k \in I_e}$ of $G_k$ bitorsors, there are $G$-bitorsors $\{P_i\}_{i \in I_v}$
  such that $P_i|_{F_k} \wedge^{G_k} P_j|_{F_k} \simeq P_k$ holds whenever $(i,j,k) \in S_I$.
\end{defn}

We will now recall a description of a $(G,H)$-bitorsor $P$ as a left $G$-torsor
together with a $G$-equivariant morphism $P \rightarrow \Isom(H,G)$.
This description was first given by Breen, cf. \cite[Section 2]{Breen1990} or  \cite[Section 1.4]{Breen2005}.
Details and proofs to the statements  made below can be found there for the case $G=H$. The proofs of the
case $G \neq H$ go through verbatim. 
Let $P$ be a $(G,H)$-bitorsor. Let $p$ be a section of $P$ over some $U \in \Sch{F}$.
Define an isomorphism of group sheaves 
\begin{align*}
  u_p \colon H|_U \rightarrow G|_U
\end{align*}
via 
\begin{align*}
  p|_V h=u_{p}(V)(h)p|_V
\end{align*}
for $V \in \Sch{U}$ and $h \in H|_U(V)$.\par
Let $p'$ be another section over the same $U$.
There is some $\gamma \in G(U)$ such that $p'=\gamma p$.
Then, we have 
\begin{align*}
  u_p'=i_{\gamma}u_p
\end{align*}
where $i_{\gamma}$ is conjugation $(g\mapsto \gamma g \gamma^{-1})$. \par
Thus, we obtain an isomorphism of group sheaves
\begin{align*}
  u \colon P &\longrightarrow \Isom(H,G) \\
  p & \mapsto u_p
\end{align*}
that is equivariant with respect to the conjugation map 
\begin{align*}
  i \colon G \longrightarrow \Aut(G).
\end{align*}
This process is reversible.
\begin{lemma}[{cf.  \cite[Lemma 2.5]{Breen1990}}] \label{right_action}
  For a left $G$-torsor $P$, the following are equivalent: 
  \begin{enumerate}
    \item the data of a right $H$ action making it a bitorsor,
    \item an  morphism of sheaves $ u \colon P \longrightarrow \Isom(H,G)$ equivariant with respect to
      $i \colon G \rightarrow \Aut(G)$.
  \end{enumerate}
\end{lemma}

The following Lemma continues this equivalence with respect to morphisms. 
\begin{lemma}[{cf.  \cite[Section 1.4.3]{Breen2005}}]
  Let $f \colon P \rightarrow P'$ be a morphism of left $G$-torsors.
  Assume that $P$ and $P'$ are $(G,H)$-bitorsors and let 
  let $u,u'$ denote the equivariant morphism of sheaves discussed in Lemma \ref{right_action}.
  Then, the following are equivalent: 
  \begin{enumerate}
    \item $f$ is a morphism of bitorsors
    \item $u=u'f$
  \end{enumerate}
\end{lemma}

 If we fix a cover  $Y \rightarrow F$ of $F$ and a section $p \in P(Y)$,
we can reinterpret the two lemmas above: 
\begin{lemma}[{cf.  \cite[Proposition 2.5]{Breen1990} and \cite[Section 1.4.5.]{Breen2005}}]\label{bitorsor}
  Let $(P,p)$ be a left $G$-torsor with a section $p$ over a cover $Y \rightarrow F$.
  Then, the following data are equivalent 
  \begin{enumerate}
    \item A right $H$-action on $P$ making it a $(G,H)$-bitorsor
    \item a sheaf isomorphism $u_p \colon H|_Y \rightarrow G|_Y$
  \end{enumerate}
  Let now $(P',p')$ be another such tuple where $p'$ is also a section over $Y$.
  Assume that both $P$ and $P'$ are $(G,H)$-bitorsors. Let $f \colon P \rightarrow P'$
  be a morphism of left $G$-torsors. Let $g \in G(Y)$ be such that 
  $f(p)=gp'$ holds. Then, the following are equivalent 
  \begin{enumerate}
    \item $f$ is a morphism of bitorsors
    \item $u_p=i_gu_p'$
  \end{enumerate}
\end{lemma}

\subsection{Patching for bitorsors} \label{sec:bit-pat}
In this section, fix a base field $F$ and let $\mathcal{C}=\Sch{F}$ equipped with the big  \'{e}tale topology.

\begin{lemma}\label{glueu}
  Let $K$ be a finite separable extension of $F$ and let $\K$ denote the corresponding
  inverse system obtained from base change. 
   Let $G,H$ be linear algebraic groups over $K$.
   Given isomorphism of group schemes
  $u_i \colon G|_{K_i} \rightarrow H|_{K_i}$ for all $i \in I_v$ such that
  \begin{align*}
    u_i|_{K_k}=u_j|_{K_k}
  \end{align*}
  whenever $(i,j,k) \in S_I$, there is an isomorphism of group schemes $u \colon G|_{K} \rightarrow H_{K}$
  satisfying $u|_{K_i}=u_i$ for all $i \in I$.
\end{lemma}
\begin{proof}
  Let $A$, $B$ be $K$-algebras such that $G|_K=\Spec(A)$ and $H|_K=\Spec(B)$.
  Then, the $u_i$ induce ring isomorphisms $f_i \colon B_i \rightarrow A_i$.
  Fix some $b \in B$. Consider the elements $f_i(b)$ for $i\in I$.
  We have $f_i(b)=f_j(b) \in A_k$ for $(i,j,k) \in S_I$ by assumption. Hence, the elements
  $f_i(b)$ determine a unique element $f(b)$ in $A$ which define a morphism $f \colon B \rightarrow A$.
  It is clear that this is a ring isomorphism. Hence, we get an isomorphism
  of schemes $u \colon G\rightarrow H$. For $u$ to be compatible with the group
  structure, we need $f$ to be compatible with Hopf algebra structure. But,
  as each $f_i$ is compatible with the Hopf algebra structure, it follows that $f$
  is as well.
\end{proof}

Let $\F$ be an inverse system of fields with limit $F$. Let $\BPAP(G,H)(\F)$ denote the category
of $(G,H)$-bitorsor patching problems over $\F$, defined analogously to $\TPAP(G)(F)$. There is a natural functor 
\begin{align*}
  \beta_{(G,H)}'' \colon \Bitorsor(G,H)(F) \rightarrow  \BPAP(G,H)(\F).
\end{align*}
We recall the natural functor
\begin{align*}
  \beta_{G}' \colon \Tor(G)(F) \rightarrow  \TPAP(G)(\F).
\end{align*}
We will now see that we can patch bitorsors whenever we can patch torsors.

\begin{theorem}
  Assume that patching holds for $G$-torsors over $\F$. Then, patching holds
  for $(G,H)$-bitorsors, i.e.
  if $\beta_G'$ is an equivalence of categories, then so is $\beta_{(G,H)}''$.
  \label{bitorsor-patching}
\end{theorem} 
\begin{proof}
  We need to prove essential surjectivity and fully faithfulness.
  Let us start with essential surjectivity. \par
  Fix some $\cP \in \BPAP(G,H)(\F)$. We have a commuting diagram of functors 
  \begin{center}
    \begin{tikzcd}
     \Bitorsor(G,H)(F)\arrow{r}{\beta_{(G,H)}''} \arrow{d} & \BPAP(G,H)(\F) \arrow{d}  \\
     \Tor(G)(F) \arrow{r}{\beta_G'}& \TPAP(G)(\F),
    \end{tikzcd}
  \end{center}
  where the vertical functors are the forgetful ones.\par
  By assumption $\beta_G'$ is essentially surjective, so there is some left $G$-torsor $P$ defined over $F$
  together with isomorphisms $\phi_i \colon P|_{F_i} \rightarrow  P_i $ for all $i \in I_v$ that are compatible
  with the morphisms $\nu_{k}$. Let $K/F$ be a finite separable field extension, such that $P(K) \neq \emptyset$.
  Fix some $p_0 \in P(K)$. Let $K_i := F_i \otimes_F K$ for all $i \in I$. 
  Set $p_i=\phi_i(K_i)(p_0|_{K_i})$, which defines a trivializing family $\{ p_i \}_{i \in I}$. Note that, by construction,
  $\nu_{k}(p_i|_{K_k})=p_j |_{K_k}$.
  By Lemma \ref{bitorsor}, we get sheaf isomorphisms $u_{p_i} \colon H|_{K_i}\rightarrow G|_{K_i}$ for all $i \in I$. Recall that this morphism is defined over $V \rightarrow K_k$ via
  $p_i|_V .h = u_{p_i}(V)(h).p_i|_V$ for $h \in H(V)$. We claim that $u_{p_i}|_{K_k}=u_{p_j}|_{K_k}$ for all $(i,j,k) \in S_I$. This follows from 
  \begin{align*}
    u_{p_i}(V)(h).p_j|_V &=u_{p_i}(h).\nu_{k}(V)(p_i|_V) \\
                        &=\nu_{k}(V)(u_{p_i}(V)(h).p_i|_V) \\
                        &=\nu_{k}(p_i|_V.h) \\
                        &=\nu_{k}(p_i|_V).h \\
                        &= p_k.h \\
                        &= u_{p_j}(V)(h).p_j|_V.      
  \end{align*}
By Lemma \ref{glueu}, we get a global  isomorphism $u_{p_0}\colon G \rightarrow H$.
  By Lemma \ref{bitorsor}, this in turn
  equips $P$ with a $(G,H)$-bitorsor structure. \par
  In order to show that $\beta''_{(G,H)}(P)$ is isomorphic to the given bitorsor patching problem $\cP$, it is enough to show that the morphisms $\phi_i$ are in fact
  morphisms of bitorsors. By Lemma \ref{bitorsor}, this is equivalent to checking that 
  \begin{align*}
    u_{p_0}|_{K_i}=u_{p_i}
  \end{align*}
  for all $i \in I_v$ (note that $\phi_i(p_0)=p_i$). But, this is clear by construction of $u_{p_0}$. Hence, $\beta_{(G,H)}''$ is essentially
  surjective.\par
  Let us now show that $\beta_{(G,H)}''$ is fully faithful. It is clearly faithful as $\beta_G'$ is faithful.
  So, we only need to prove that it is full. Fix two $(G,H)$-bitorsors $P,P'$ over $F$ and let $\beta_{(G,H)}''(P)=\cP$ and $\beta_{(G,H)}''(P')=\cP'$.
  Let $\alpha \colon \cP \rightarrow \cP'$ be a morphism
  in $\BPAP(G,H)(\F)$. Note that $\alpha$ is also a morphism in $\TPAP(G)(\F)$. Hence, there
  is a morphism $\alpha \colon P \rightarrow P'$ of left $G$-torsors inducing the morphisms
  $\alpha_i \colon P_i \rightarrow P_i'$. A straightforward check shows that $\alpha$ is a morphism of bitorsors.
\end{proof}

In the context of gerbe patching, we will have to patch bitorsors over covers of $F$. These covers are formed
by finite products of finite separable field extensions. We thus need to extend our patching results
to this setup.

\begin{corollary}\label{product_patching}
Let $A$ be a finite product of finite separable field extensions of $F$ and let
$G,H$ be linear algebraic groups over $A$.
If vector space patching holds over $\F$ then $(G,H)$-bitorsor patching holds over  $\F_A$.
\end{corollary} 
\begin{proof}
  Follows from Theorem \ref{bitorsor-patching} and
  Proposition \ref{trs}.
  (Note that the proof of Theorem \ref{bitorsor-patching} goes through verbatim
  if we replace $\F$ by $\F_A$.)
\end{proof}

\subsection{A Mayer-Vietoris sequence and a local-global principle for bitorsors}

Using patching for bitorsors, we can construct the start of our Mayer-Vietoris sequence, which will be extended
in section \ref{gerbe-lg}. Fix a group sheaf $G$ and let $Z$ denote its center.

We say that $G$-bitorsors satisfy \emph{local-global principle over $\F$} if for any $G$-bitorsor $T$,
we have  that $T \simeq G$ if and only if $T_i \simeq G_i$ for all $i \in I_v$.

Recall that automorphisms of the trivial $G$-bitorsor $G$ can be identified with $Z(F)$. Furthermore,
note that
 $\HH^{-1}(F,G \rightarrow \Aut(G))=Z(F)$ and that $\HH^0(F,G \rightarrow \Aut(G))$
classifies bitorsors up to isomorphism (cf. Appendix \ref{hyper}).

Note that we have maps
\begin{align*}
    \HH^0(F_i, G \rightarrow \Aut(G))  \times \HH^0(F_j, G \rightarrow \Aut(G)) \rightarrow
    \prod_{i \in I_k} \HH^0(F_k, G \rightarrow \Aut(G))
\end{align*}
for $(i,j,k) \in S_I$ defined via $(P_i,P_j) \mapsto P_i\wedge^G P_j^{\op}$.
This induces a map
\begin{align*}
\prod_{i \in I_v} \HH^0(F_i, G \rightarrow \Aut(G)) \rightarrow
    \prod_{i \in I_k} \HH^0(F_k, G \rightarrow \Aut(G)).
\end{align*}

\begin{lemma}
  There is a map of pointed sets
  \begin{align*}
    \prod_{i \in I_k} \HH^{-1}(F_k, G \rightarrow \Aut(G)) \rightarrow
      \HH^0(F, G \rightarrow \Aut(G))
  \end{align*}
\end{lemma} 
\begin{proof}
   We can define the map 
as follows: Given elements $e_k \in Z(F_k) $, consider the $G$-bitorsor patching
problem $(\{ G_i\}_i, \{ e_k \}_k )$. By Theorem \ref{bitorsor-patching}, there is a $G$-bitorsor
$P$ over $F$ in the essential preimage of the patching problem.
Map  $(e_k)_k$ to the equivalence class of $P$. Note that this
is independent of the choice of $P$.
\end{proof}

\begin{theorem}[Mayer-Vietoris for non-abelian hypercohomology (1)] \label{mayer1}
  Assume that $G$-bitorsor patching holds over $\F$. Then, there is an exact sequence
  \begin{adjustwidth}{-2em}{-2em}
  \begin{center}
    \begin{tikzcd}[row sep =small, column sep = small]
  1 \arrow{r} &    \HH^{-1}(F, G \rightarrow \Aut(G)) \arrow{r}&
    \prod_{i \in I_v} \HH^{-1}(F_i, G \rightarrow \Aut(G))\arrow{r} \arrow[draw=none]{d}[name=Z, shape=coordinate]{} & 
    \prod_{k \in I_e} \HH^{-1}(F_k, G \rightarrow \Aut(G))  \arrow[rounded corners,to path={ -- ([xshift=2ex]\tikztostart.east)|- (Z) [near end]\tikztonodes-| ([xshift=-2ex]\tikztotarget.west)-- (\tikztotarget)}]{dll} \\
 \; &     \HH^0(F, G \rightarrow \Aut(G)) \arrow{r}& 
    \prod_{i \in I_v} \HH^0(F_i, G \rightarrow \Aut(G))\arrow{r} \arrow{r}&
    \prod_{k \in I_e} \HH^0(F_k, G \rightarrow \Aut(G))
    \end{tikzcd}
  \end{center}
  \end{adjustwidth}
\end{theorem} 
\begin{proof}
  Exactness in the first row follows from the description $\HH^{-1}(F,G \rightarrow \Aut(G))=Z(F)$ and the fact that
  $\F$ is a factorization inverse system with limit $F$.  Exactness at $\prod_{i \in I_v} \HH^0(F_i, G \rightarrow \Aut(G))$
  follows from Theorem \ref{bitorsor-patching}. Finally, exactness at $\HH^0(F, G \rightarrow \Aut(G))$ also follows from
  Theorem \ref{bitorsor-patching} and the fact that the automorphism group of $G$ as a bitorsor is $Z$.
\end{proof}

The sequence above allows us to deduce a local-global principle for bitorsors with respect to patches.

\begin{corollary} \label{local-global-bitorsor}
  Assume that $G$-bitorsor patching holds over $\F$.
  Then, $G$-bitorsors satisfy local-global principle over $\F$ iff
  $Z(G)$ satisfies factorization over $\F$.
\end{corollary}

\section{Gerbes}
The goal of this section is to prove our main theorem on gerbe patching.

\subsection{Preliminaries}
Let us again fix the big \'{e}tale site $\SchF$ of the scheme $\Spec(F)$ for some fixed field $F$.
As before, we will often just write $F$ to denote $\Spec(F)$.
Note that covers of $F$ can all be taken to be of the form $Y \rightarrow F$ for a single arrow where $Y$ is a finite product
of finite separable field extensions of $F$. We will use this
throughout this section as it will simplify notation.

We will first recall the definition of a gerbe as a special kind of stack. Our main reference for
 gerbes is \cite{Giraud1971}.

\begin{defn}
  Let $\G \rightarrow \SchF$ be a stack. We say that $\G$ is a \textit{gerbe} if  the following two conditions are satisfied: 
  \begin{enumerate}
    \item For all $X \in \SchF$, there is a cover $Y \rightarrow X$ such that $\G(Y) \neq \emptyset$
    \item For any $x,x' \in \G(X)$, there is a cover $Y \rightarrow X$ such that the pullbacks of
      $x$ and $x'$ are isomorphic in $\G(Y)$.
  \end{enumerate}
\end{defn}
 
\begin{example}
  Let $G$ be a sheaf of groups over $\SchF$. Then, the classifying stack $BG$ is a gerbe.
  Its fiber over $Z \rightarrow F$ is given by $BG(Z)=\Tor(G)(Z)$.
\end{example}

\begin{defn}
  We say that a gerbe $\G \rightarrow \SchF$ is \textit{trivial} (or \textit{neutral}) if
  $\G(F) \neq \emptyset$.
\end{defn}
\begin{example}
  The gerbe $BG$ is trivial as the trivial left $G$-torsor $G$ is in $BG(F)$.
\end{example}

\begin{remark}
  Let $x,y \in \G(Z)$ for some $Z \in \SchF$. The sheaf
  $\Isom_Z(x,y)$ on $\Sch{Z}_{\acute{e}t}$ is a 
  $(\Aut_Z(y),\Aut_Z(x))$-bitorsor.
\end{remark}

Recall that  a 2-category $\calC$ consists of a class of objects $\obj(\calC)$ together with
categories $\Hom(A,B)$ for $A,B \in \obj(\calC)$ and composition functors
$\Hom(A,B) \times \Hom(B,C) \rightarrow \Hom(A,C)$ satisfying the usual identities.
Additionally, we require $\Hom(A,A)$ to contain an identity element with respect
to composition. Similarly, for $f \in \Hom(A,B)$, we require there to be an identity
element in $\Hom(f,f)$.
We will call an object in $\Hom(A,B)$ a \emph{1-morphism} and a morphism in $\Hom(A,B)$ a \emph{2-morphism}.

 A typical example of a 2-category is the category of categories,
where for two categories $\mathcal{A},\mathcal{B}$, the category $\Hom(\mathcal{A},\mathcal{B})$ has functors
as objects and natural transformations as morphisms.

As stacks are in particular categories, we consider the 2-category of stacks over some site as the category with
objects being stacks, 1-morphisms being functors of fibered categories and 2-morphisms being natural transformations
of functors.

Similarly, we will consider the 2-category of gerbes in this section and will later define the 2-category of
gerbe patching problems.

\begin{defn}
  A 1-morphism of gerbes $\G \rightarrow \G'$ over $\SchF$ is a morphism
  of stacks $\G \rightarrow \G'$ over $\SchF$. \par
  A 2-morphism of gerbes is a natural transformation of functors. \par
  By $\Gerbe(F)$, we denote the 2-category of gerbes with 1-morphisms given by
 equivalences and 2-morphisms given by natural isomorphisms of
  functors.
\end{defn}

The category $\Hom(BG,BH)$ in $\Gerbe(F)$ is equivalent to the
category of $(G,H)$-bitorsors.

\begin{theorem}[{\cite[Giraud]{Giraud1971}}] \label{morita}
\label{torsorequi}
  Let $P$ be a $(H,G)$-bitorsor over $F$. Then, the functor 
  \begin{align*}
    BG& \rightarrow BH \\
    (S,T) &\mapsto (S, P_S \wedge^{G_S} T)
  \end{align*}
  is an equivalence of categories. Furthermore, any equivalence between
  $BG$ and $BH$ is of this form.
\end{theorem} 
\begin{remark}
  Given two equivalences $BG \rightarrow BH$ and $BH \rightarrow BK$ for
 group sheaves $G,H$ and $F$,
  we obtain an $(H,G)$-bitorsor $P$ and an $(K,H)$ bitorsor $P'$. 
The composition $BG \rightarrow BH \rightarrow BK$ then corresponds
  to $P' \wedge^H P$.
\end{remark}

\subsubsection{A semi-cocyclic description of a gerbe}

We will now follow sections 2.3-2.6 of \cite{Breen1991} to introduce a cocyclic description of a gerbe $\G \rightarrow \SchF$.

Let $Y \rightarrow F$ be a cover such that there is an element $y \in \G(Y)$. Let $G=\Aut_Y(y)$ denote
the sheaf of automorphisms of $y$ over $\Sch{Y}$.  The choice of an object $y$ defines an equivalence 
\begin{align*}
   \Phi\colon \G|_Y \longrightarrow BG 
\end{align*}
defined on fibers over $f \colon Z \rightarrow Y$ via 
\begin{align*}
  \Phi(Z) \colon \G(Z) & \longrightarrow \Tor(G)(Z) \\
  z & \mapsto \Isom_Z(f^*y,z)                   
\end{align*}
where $\Isom_Z(z,f^*y)$ is a $G_y|_Z$-torsor by the natural action on the left. \par
We fix the following notation: Let $\pr_i \colon Y^2 =Y \times_F Y \rightarrow Y$ denote the natural
projections on the $i$-th factor for $i=1,2$. We will denote by $G_i=\pr_i^*G$ the pullbacks of $G$
to $Y^2$.
Let $\pr_{ij} \colon Y^3\rightarrow Y^2$ denote the natural projections on
the $i$-th and $j$-th component for $1\leq i <j \leq 3$.  We will denote by $G^i$ the pullbacks of
$G$ to $Y^3$ along the natural projections $Y^3 \rightarrow Y$.
Finally, let $\pr_{ijk} \colon Y^4 \rightarrow Y^3$
denote the natural projection on the $i$-th, $j$-th and $k$-th component for $1\leq i <j <k \leq 4$. 
Let $G_{(i)}$ denote the pullbacks of $G$ to $Y^4$ along the natural projections $Y^4 \rightarrow Y$.
Note that $\Phi$ induces an equivalence
\begin{align*}
\varphi=\pr_1^* \Phi \circ \left(\pr_2^* \Phi\right)^{-1} \colon BG_2 \longrightarrow BG_1. 
\end{align*}
where $G_i=\pr_i^* G$. \par
We obtain an isomorphism of functors
 \begin{align*}
   \psi \colon \pr_{12}^*\varphi \circ \pr_{23}^*\varphi \Rightarrow \pr_{13}^*\varphi                                      
 \end{align*}
rather than an equality. The various pullbacks to $Y^4$ of this isomorphism make the following diagram of natural
transformations commute:
 \begin{center}
  \begin{tikzcd}[column sep=large, row sep=large]
  BG|_{Y^4} \arrow{rrr}[name=F, below]{} \arrow{rd}[name=T]{} \arrow{rrd}[name=S,below]{} 
  \arrow[Rightarrow,to path=(S) -- (F)\tikztonodes,shorten <= 0.3em]{u}{}
  \arrow[Rightarrow,to path=(T) -- (S)\tikztonodes]{u}{}
  &\; &\; & BG|_{Y^4}  \\
  \;  & BG|_{Y^4} \arrow{r}[name=Fo]{} \arrow{rru}[name=Fi,below]{}
  \arrow[Rightarrow,to path=(Fi) -- (F)\tikztonodes,shorten <= 0.3em]{u}{}& BG|_{Y^4} \arrow{ru}[name=Si]{}
  \arrow[Rightarrow,to path=(Si) -- (Fi)\tikztonodes]{u}{}& \; 
  \end{tikzcd}
\end{center}
where $Y^4=Y\times_F Y \times_F Y \times_F Y$.

The gerbe $\G$ is completely determined by the tuple $(G,\varphi,\psi)$.

We will now reinterpret the above description in terms of bitorsors.
Note that the equivalence $\varphi$ can be identified with the $(G_2,G_1)$-bitorsor $E=\Isom(\pr_2^* y, \pr_1^* y)$. 
Then, the isomorphism $\psi$ corresponds to an isomorphism of $(G^1,G^3)$-bitorsors 
 \begin{align*}
   \pr_{12}^*E \wedge^{G^2} \pr_{23}^* E \longrightarrow \pr_{13}^* E
 \end{align*}
where $G^i$ denotes the pullback of $G$ along $Y \times_F Y \times_F Y \rightarrow Y$ on the $i$-th component. \par
Finally, the compatibility condition translates to the commutativity of the diagram 
\begin{center}
  \begin{tikzcd}
   \pr_{12}^* E \wedge^{G_{(2)}} \pr_{23}^* E \wedge^{G_{(3)}} \pr_{34}^* E \arrow{r} \arrow{d}& \pr_{13}^* E \wedge^{G_{(3)}} \pr_{34}^* E \arrow{d}  \\
   \pr_{12}^* E \wedge^{G_{(2)}} \pr_{24}^* E \arrow{r} & \pr_{14}^* E 
  \end{tikzcd}
\end{center}
where the arrows are induced by the various pullbacks of $\psi$ and $G_{(i)}$ denotes the pullback of $G$ along $Y^4 \rightarrow Y$ onto the $i$-th component.
So, a gerbe can be described by the triple $(G,E,\psi)$ and we will henceforth go back and forth between the categorical and the cocyclic description. We call the triple $(G,E,\psi)$ the \emph{semi-cocyclic description of $\G$}.

Let us now turn to morphisms of gerbes.

Let $\rho \colon \G \rightarrow \G'$ be a morphism of gerbes over a field $F$. Let $K/F$ be a finite separable extension such that there are
$x \in \G(K)$ and $y \in \G'(K)$. After possibly replacing $K$ by a further finite separable extension, we may assume that
$\rho(x)$ and $y$ are isomorphic in $\G'(K)$.
Using the cocyclic description, we get $\G=(G,E,\psi)$ and $\G'=(G',E',\psi')$ where we choose the descriptions
induced by $x$ and $y$. In particular, $\Aut(x)=G$ and $\Aut(y)=G'$.

Note that $\rho$ induces a map 
\begin{align*}
   \rho' \colon G \longrightarrow G'                                 
\end{align*}
and a $\rho'|_{K \times_F K}$-equivariant  map 
\begin{align*}
  \alpha \colon E \rightarrow E'.
\end{align*}

Note that $\alpha$ is compatible with $\psi$ and $\psi'$. 
Conversely, given a map $\alpha$ compatible with $\psi$ and $\psi'$, one gets an induced morphism
of gerbes $\G \rightarrow \G'$.

\begin{lemma}[{\cite[Section 2.6]{Breen1991}}]\label{gerbe_morphism}
  Let $\rho \colon \G \rightarrow \G'$ be a morphism of gerbes. Assume that there is
  a finite separable cover $K/F$ such that $\G(K)\neq \G'(K)$. Fix some
  $x \in \G(K)$ and $y \in \G(K)$. Assume that there is a morphism
  $\rho(K)(x) \rightarrow y $ in $\G'(K)$.
  Rewrite $\G=(\Aut(x),E,\psi)$ and $\G=(\Aut(y),E',\psi')$ using the construction
  described above coming from $x$ and $y$. Then, $\rho$ induces an equivariant isomorphism
  of bitorsors $E \rightarrow E'$ that is compatible with $\psi$ and $\psi'$.
  Conversely, given an isomorphism $E \rightarrow E'$ compatible with the gluing data, one can construct
  a morphism of gerbes.
\end{lemma}

\subsection{Bands and patching of non-abelian $H^2$} \label{bands}

In this section, we will continue to work over the big \'{e}tale site of schemes over $F$.
Before we can patch gerbes, we will first investigate when we can patch equivalence classes of gerbes, i.e.
elements in the non-abelian second cohomology set of a band.

Given a gerbe $\G$, one can associate to it a band $L$. Let us quickly  review the definition of a band.
For more details, we refer again to \cite{Giraud1971}, compare also the appendix of \cite{Deligne1982}.

Let $H,G$ be group sheaves and  $G^{\ad}=G/Z$, where $Z$ is the center of $G$. Then, $G^{\ad}$ acts on
$\Isom(H,G)$ via conjugation. Let $\Isex(H,G)$ denote the quotient of $\Isom(H,G)$ by the $G^{\ad}$-action.
Then, a \emph{band} $L$ consists of a triple $(Y,G,\varphi)$ where $Y$ is a cover of $F$, $G$ is a group sheaf
defined over $Y$ and $\varphi \in \Isex(G_1,G_2)$ where $G_i= \pr_i G$ for $\pr_i \colon Y \times_F Y \rightarrow Y$.
We require $\varphi$ to satisfy the cocycle condition $\pr_{31}^* \varphi = \pr_{32}^* \varphi \circ \pr_{21}^* \varphi$ over $Y\times_F \times_F Y$.
Thus, a band is a descent datum for a group sheaf over a cover $Y \rightarrow F$  modulo inner automorphisms. Every group sheaf over $F$ defines a band with trivial
gluing datum. Furthermore, any abelian band (i.e. $G$ is abelian) is induced by a group sheaf over $F$ as the datum of a band in this case just gives descent
datum since inner automorphisms of $G$ are all trivial.
The center $Z(L)$ of band $L=(Y,G,\varphi)$ is defined as $(Y,Z,\varphi_{Z})$ and we identify it with the group sheaf
over $F$ determined by the descent datum of the band.
For a cover $f\colon Y' \rightarrow Y$, we identify the bands $(Y,G,\varphi)$ and $(Y',f^*G,f^*\varphi)$.
An isomorphism between two bands $(Y,H,\varphi) \rightarrow (Y,G,\tau)$ is given by an element $g \in \Isex(H,G)$ compatible
with $\varphi$ and $\tau$.

We say that a band $L=(Y,G,\varphi)$ is \emph{linear algebraic}, if $G$ is a linear algebraic group, i.e.
if $Y=\Spec(R)$ for $R=\prod_i L_i$ where $L_i$ is a finite separable field extension, $G_{L_i}$ is a linear algebraic group.

Given a gerbe $\G$, we can define an associated band. Pick a cover $Y \rightarrow F$ and an object $x \in \G(Y)$. Let $G = \Aut(x)$ be the sheaf of automorphisms of $x$. Let $x_i = \pr_i^* x$ denote the two pullbacks along the projections 
$\pr_i \colon Y \times_F Y \rightarrow Y$. By the definition of gerbes, there is a cover $U \rightarrow Y\times_F Y$
such that $x_1|_U$ and $x_2|_U$ are isomorphic in $\G(U)$. An isomorphism $ f \colon x_1|_U \rightarrow x_2|_U$ defines an isomorphism
$\lambda_f \colon G_1|_U \rightarrow G_2|_U$. If $g \colon x_1|_U \rightarrow x_2|_U$ is another isomorphism,
then $\lambda_f,\lambda_g$ differ by an inner automorphism of $G_2$. Thus, there is a well defined element $\lambda \in \Isex(G_1,G_2)(U)$.
An easy calculation shows that the pullbacks of $\lambda_f$ on $U\times_{Y\times_F Y} U$ differ by an inner automorphism. 
Hence, the pullbacks of $\lambda$ agree on $U\times_{Y\times_F Y} U$ and we thus obtain an element $\lambda' \in \Isex(G_1,G_2)(Y \times_k Y)$ whose restriction to $U$ equals $\lambda$.
It is not hard to see that $\lambda'$ satisfies the cocycle condition and thus $(Y,G,\lambda')$ defines
a band. Given a gerbe, the associated band is unique up to unique isomorphism and we denote it by $\Band(\G)$.
Note that a morphism $\G \rightarrow \G'$ induces a map $\Band(\G) \rightarrow \Band(\G')$.
Given a band $L$, a \emph{gerbe banded by $L$} is a tuple $(\G,\theta)$ where $\G$ is a gerbe
and $\theta \colon \Band(\G) \rightarrow L$ is an isomorphism of bands. We often suppress $\theta$ from the
notation.

Given two gerbes $\G,\G'$ banded by $L$, an \emph{$L$-morphism of gerbes} is a morphism $\alpha \colon\G \rightarrow \G'$ that is compatible with
the morphisms $\theta$ and $\theta'$, i.e. the diagram
\begin{center}
  \begin{tikzcd}
\Band(\G) \arrow{rr}{\alpha}   \arrow{dr}{\theta}  &  & \Band(\G') \arrow{dl}{\theta'} \\
    & L &
  \end{tikzcd}
\end{center}
commutes. Any morphism of $L$-gerbes is an equivalence.

Let $\HH^2(F,L)$ denote the set of $L$-equivalence classes of gerbes banded by $L$ (this means that we only consider equivalences
of $L$ gerbes that are compatible with the band $L$). If $L$ is abelian coming from
the group sheaf $A$ over $F$, then this definition coincides with the usual definition of $\HH^2(F,A)$ 
as a Galois cohomology group.

Furthermore, there is a remarkable relation between $\HH^2(F,L)$ and $\HH^2(F,Z(L))$.
\begin{theorem}[{\cite[Theorem 3.3.3]{Giraud1971}}] \label{action}
  If $\HH^2(F,L)$ is not empty, then it is a principal homogeneous space under $\HH^2(F,Z(L))$.
\end{theorem}

Harbater, Hartmann and Krashen have proved patching for Galois cohomology groups of abelian algebraic groups over arithmetic curves under some
mild compatibility assumptions between the characteristic of $F$ and the order of the group (compare \cite{Harbater2012a}).
Using the above theorem, we can deduce patching for non-abelian Galois cohomology over arithmetic curves whenever the center
of the band admits patching. We will pursue this in Section \ref{curve}. For the purpose of gerbe patching, let us note this
immediate consequence.

We define patching for $\HH^2(\circ,L)$ over a finite inverse factorization system $\F$ analogously to the case of patching
for Galois cohomology, see Section \ref{gal}. Note that for $(i,j,k) \in S_I$, we have two maps
$\HH^2(F_i,L) \times \HH^2(F_j,L) \rightrightarrows \HH^2(F_k,L)$ given by restriction of the first and the 
second factor respectively. Hence, we have two maps $\prod_{i \in I_v} \HH^2(F_i,L)  \rightrightarrows \prod_{k \in I_e} \HH^2(F_k,L)$.

\begin{defn}
  Let $\F$ be a finite inverse factorization system with inverse limit $F$. For a band $L$ over $F$, we say that
\emph{patching holds for $\HH^2(\circ,L)$ over $\F$} if the the following sequence is an equalizer diagram
\begin{align*}
  \HH^2(F,L) \rightarrow \prod_{i \in I_v} \HH^2(F_i,L)  \rightrightarrows \prod_{k \in I_e} \HH^2(F_k,L).
\end{align*}
\end{defn}

\begin{proposition} \label{patching-coho}
  Let $\F$ be a finite inverse system for fields with inverse limit $F$ and let $L$ be a band over $F$.
  If patching holds for $\HH^2(\circ,Z(L))$ over $\F$ and $\HH^2(F,L) \neq \emptyset$, then patching holds for
  $\HH^2(\circ,L)$ over $\F$.
\end{proposition} 
\begin{proof}
  By assumption, there is a class $\alpha \in \HH^2(F,L)$. Given a patching problem $\{\beta_i\}_{i \in I_v}$ with
  $\beta_i \in \HH^2(F_i,L)$, there are elements $\gamma_i \in \HH^2(F_i,Z(L))$ such that $\beta_i=\gamma_i. \alpha|_{F_i}$
  by Theorem \ref{action}. As the action of $\HH^2(F_i,Z(L))$ on $\HH^2(F_i,L)$ is simply transitive,
  it follows that $\{\gamma_i \}_{i \in I_v}$ defines a patching problem for $\HH^2(\circ,Z(L))$. By assumption,
  there is $\gamma \in \HH^2(F,Z(L))$ such that $\gamma|_{F_i} =\gamma_i$ for all $i \in I_v$. The element
  $\beta := \gamma. \alpha$ solves the patching problem $\{\beta_i \}_{i \in I_v}$.
\end{proof}

\subsection{Patching gerbes}

Let $\F$ be a factorization inverse system of fields with inverse limit $F$. Let $L$ be a band over $F$. Let $\Gerbe(F,L)$ denote the 2-category of $L$-banded gerbes $\G$ over $F$. Here, morphisms are given by equivalences of gerbes and 2-morphisms are given by natural
isomorphisms. Let $\GPAP(\F,L)$ denote the 2-category of $L$-gerbe patching problems, i.e. an object consists of  a collection of gerbes $L$-banded gerbes  $\G_i$ over $F_i$ together with $L$-equivalences
 $\sigma_{k} \colon \G_i|_{F_k} \rightarrow \G_j|_{F_k}$ for $(i,j,k) \in S_I$.

A 1-Morphisms $(\{\G_i\},\{\sigma_k \}) \xrightarrow{(\alpha,f)}(\{\G'_i\},\{\sigma'_k \}) $ is given by a collection of equivalences of gerbes
$\alpha_i \colon \G_i \rightarrow \G_i'$ and natural isomorphisms $f_{k} \colon \sigma_k' \circ \alpha_i|_{F_k} \Rightarrow \alpha_j|_{F_k} \circ \sigma_k$, 
pictorially:

\begin{center}
  \begin{tikzcd}
    \G_i|_{F_{k}}  \arrow{r}[name=U]{\alpha_i|_{F_k}} \arrow{d}[left]{\sigma_k} & \G_{i}'|_{F_k} \arrow{d}{\sigma_k'} \\
    \G_j|_{F_k}    \arrow{r}[below,name=L]{\alpha_j|_{F_k}}  & \G_{j}'|_{F_k}  \arrow[Rightarrow,to path=(U) -- (L)\tikztonodes,shorten <= 0.3em,shorten >= 0.3em]{u}{f_{k}}
  \end{tikzcd}
\end{center}
Composition of morphisms is given by composing equivalences of gerbes and by horizontally composing the natural isomorphisms.

Given two morphisms $(\alpha,f),(\beta,g) \colon \{\G_i\},\{\sigma_k \}) \rightarrow (\{\G'_i\},\{\sigma'_k \})$,
a 2-morphism $\psi=(\{\psi_i\})$ is given by a collection of natural isomorphisms

\begin{center}
  \begin{tikzcd}[column sep = huge]
   \G_i \arrow[bend left]{r}[above,name = U]{\alpha_i} \arrow[bend right]{r}[below,name = L]{\beta_i}    &  \G_i'  
   \arrow[Rightarrow,to path=(U) -- (L)\tikztonodes,shorten <= 0.3em,shorten >= 0.3em]{u}{\psi_i}
  \end{tikzcd}
\end{center}

such that the diagram of natural transformations commutes:

\begin{center}
  \begin{tikzcd}
    \G_i|_{F_{k}}  \arrow{r}[name=U]{\alpha_i} \arrow{d}[left]{\sigma_k}   \arrow[Rightarrow,shorten <= 1.5em,shorten >= 1.5em,bend left]{rrr}{\psi_i|_{F_k}} & \G_{i}'|_{F_k} \arrow{d}{\sigma_k'} 
    &  \G_i|_{F_{k}}  \arrow{r}[name=U2]{\beta_i} \arrow{d}[left]{\sigma_k}  & \G_{i}'|_{F_k} \arrow{d}{\sigma_k'}  \\
    \G_j|_{F_k}    \arrow{r}[below,name=L]{\alpha_j} \arrow[Rightarrow,shorten <= 1.5em,shorten >= 1.5em,bend right]{rrr}{\psi_j|_{F_k}} 
    & \G_{j}'|_{F_k}  \arrow[Rightarrow,to path=(U) -- (L)\tikztonodes,shorten <= 0.3em,shorten >= 0.3em]{u}{f_{k}}
    & \G_j|_{F_k}    \arrow{r}[below,name=L2]{\beta_j}  & \G_{j}'|_{F_k}  \arrow[Rightarrow,to path=(U2) -- (L2)\tikztonodes,shorten <= 0.3em,shorten >= 0.3em]{u}{g_{k}}
   
  \end{tikzcd}
\end{center}

i.e. 
\begin{center}
  \begin{tikzcd}
    \sigma_k'\circ \alpha_i|_{F_K}   \arrow[Rightarrow]{r}{f_{k}}  \arrow[Rightarrow]{d}{\psi_i'} & \alpha_j|_{F_k}\circ\sigma_k  \arrow[Rightarrow]{d}{\psi_{j}'} \\
    \sigma_k'\circ\beta_i|_{F_k}    \arrow[Rightarrow]{r}{g_{k}}  & \beta_j|_{F_k} \circ \sigma_{k}
  \end{tikzcd}
\end{center}
commutes.
Here $\psi_i'$ is the natural isomorphism induced by $\psi_i$ and $\sigma_k'$.
Composition of 2-morphisms is given by vertical composition of the various natural transformations.

\begin{lemma}\label{quasi-inverse}
  Every 1-morphism in $\GPAP(\F,L)$ admits a quasi-inverse 
\end{lemma} 
\begin{proof}
  Let  $(\alpha, f)\colon (\{\G_i\},\{\sigma_k \}) \rightarrow (\{\G'_i\},\{\sigma'_k \}) $ be a 1-morphism.
  Let $\beta_i$ denote a quasi-inverse of $\alpha_i$. Then, we get natural isomorphisms $g_k\colon \beta_j|_{F_k} \circ \sigma'_k \Rightarrow \sigma_k \circ \beta_i|_{F_k}$ defined via
  \begin{align*}
    \beta_j|_{F_k} \circ \sigma'_k \Rightarrow \beta_j|_{F_k} \circ \sigma'_k \circ \alpha_i|_{F_k}\circ \beta_{i}|_{F_k} \Rightarrow \beta_j|_{F_k} \circ \alpha_j|_{F_k} \circ \sigma_k \circ \beta_{i}|_{F_k}
    \Rightarrow \sigma_k \circ  \beta_{i}|_{F_k}.
  \end{align*}
  Here, the first and last natural transformation are induced by  fixed natural isomorphisms $\phi_i \colon \alpha_i \circ \beta_i \Rightarrow \id $, while the middle arrow is induced by $f_k$.
  Hence, we can define a 1-morphism  $(\beta, g)\colon (\{\G'_i\},\{\sigma'_k \}) \rightarrow (\{\G_i\},\{\sigma_k \})$.
  It remains to check that $(\alpha,f) \circ (\beta,g)$ and $(\beta,g) \circ (\alpha,f)$ are 2-isomorphic to the identity morphism.
  We will prove that $(\alpha,f) \circ (\beta,g)$ is isomorphic to the identity, the other case is analogous.
  For this, we need to give a collection of 2-isomorphisms $\psi_i \colon \alpha_i \circ \beta_i \Rightarrow \id$ that are compatible with
  $g_k$ and $f_k$. It is tedious but straightforward to check that the choice $\psi_i=\phi_i$ works.
\end{proof}

Note that there is a natural functor of 2-categories
\begin{align*}
  \beta_{L}''' \colon \Gerbe(F,L) \longrightarrow  \GPAP(\F,L)
\end{align*}
induced by base change. We say that \emph{patching holds for $L$-gerbes over $\F$} if
$\beta_{L}'''$ is an equivalence.

In order to prove our main result, we need the following elementary lemma from the general theory of stacks.
Given two stacks $\XX$ and $\YY$ be stacks  in groupoids over $\Sch{F}$,
let $\{ U \rightarrow X\}$ be a cover in $\C$. \par
Let $\pr_i \colon U \times_X U \rightarrow U$ and $\pr_{ij} \colon U\times_X U \times_X U \rightarrow U\times_X U$ denote
the usual projections.
\begin{lemma} \label{stacks}
  The following data are equivalent: 
  \begin{enumerate}
    \item a morphism of stacks $\XX \rightarrow \YY$
    \item a morphism of stacks $\alpha \colon \XX|_U \rightarrow \YY|_U$
      together with a natural transformation $\psi \colon \pr_1^* \alpha \rightarrow \pr_2^* \alpha$
      such that the diagram
      \begin{center}
        \begin{tikzcd}
          \pr_{12}^*\pr_1^* \alpha  \arrow{r}{\pr_{12}^*\psi} \arrow[equal]{d}  & \pr_{12}^*\pr_2^* \alpha  \arrow[equal]{r} & \pr_{23}^*\pr_1^* \alpha \arrow{d}{\pr_{23}^*\psi} \\
          \pr_{13}^*\pr_1^* \alpha \arrow{r}{\pr_{13}^*\psi}   & \pr_{13}^*\pr_2^* \alpha \arrow[equal]{r}  & \pr_{23}^*\pr_2^* \alpha
        \end{tikzcd}
      \end{center}
      commutes.
  \end{enumerate}
\end{lemma}  
\begin{proof}
  This immediately follows from the fact that the category of morphisms $\XX \rightarrow \YY$ is a stack.
\end{proof}

\begin{defn}
  Let $(\{\G_i\},\{\nu_k \})$ be a gerbe patching problem. We say that
   $(\{\G_i\},\{\nu_k \})$ \emph{has property D} if  there is a cover $Z \rightarrow F$ such that
      $\G_i(Z_i) \neq \emptyset$ for all $i \in I_v$ and that there are elements $x_i \in \G_i(Z_i)$ such that
      $\nu_k(x_i|_{Z_k})$ is isomorphic to $x_j|_{Z_{k}}$ for all $(i,j,k) \in S_I$.
\end{defn}

\begin{proposition} \label{object}
  Let $L=(Y,G,\psi)$ for some cover $Y \rightarrow F$. Let $(\{\G_i\},\{\nu_k \})$ be an $L$-gerbe patching problem
  with property D. Then, if patching holds for $G$-torsors, there is $\G \in \Gerbe(F,L)$ such that $\beta_{L}'''(\G) \simeq (\{\G_i\},\{\nu_k \})$.
\end{proposition}
\begin{proof}
Since the gerbe patching problem has property D, there is  a cover $Z \rightarrow F$ and elements
$x_i \in \G_i(Z_i)$ such that $\nu_k(x_i|_{Z_k}) \simeq x_j|_{Z_k}$ for $(i,j,k) \in S_I$. We may assume without loss of generality that 
$Z \rightarrow F$ factors through $Y \rightarrow F$ as we could
replace $Z$ by $Z \times_F Y$. Hence, we may assume $L=(Z,G,\psi)$ and that we can patch $G$-torsors over $Z$.
  Let $\pr_j \colon Z \times_{F} Z \rightarrow Z$ denote the natural projection
  for $j=1,2$. Let $G_j= \pr_j^{*} G$. \par
  Let  $A = Z\times_F Z$ 
  and note that patching holds for $(G_1,G_2)$-bitorsors  over $\F_A$ by assumption
  and Corollary \ref{product_patching}. \par
  By use of the cocyclic description of gerbes and their morphisms, we will show that the given gerbe patching problem induces a bitorsor patching
  problem in $\BPAP(G_2,G_1)(\F_A)$.
Describing the gerbes $\G_i$  with respect to $Z_i$ as a cocycle,
we get tuples $(G|_{Z_i},P_i,\psi_i)$ where $P_i$ are $(G_1,G_2)$-bitorsors over
$A_i$.
 The equivalences of  the various $\G_i|_k$ with the  $\G_j|_k$ for $(i,j,k) \in S_I$ translate to
isomorphisms $P_i|_{A_k} \rightarrow P_j|_{A_k}$ by Lemma \ref{gerbe_morphism}.
 Hence, we get an element in
in $\BPAP(G_2,G_1)(\F_A).$

Thus, there is a $(G_2,G_1)$-bitorsor $P$ defined over $A=Z \times_F Z$. 
The morphisms $\psi_i$ from the cocyclic description of $G_i$ glue together by
bitorsor patching to give a global isomorphism $\psi\colon P_{12} \wedge^{G^2} P_{23} \rightarrow P_{13}$ of $(G^1,G^3)$-bitorsors. Also, again by bitorsor patching,
the morphism $\psi$ satisfies the coherence condition. Hence, we get a cocycle
$(G,P,\psi)$ defining an $L$-gerbe $\G$ in $\Gerbe(F,L)$. By construction, we obtain a 1-morphism $\beta_{K,G}'''(\G) \rightarrow (\{\G_i\},\{\nu_k\})$.
\end{proof}

\begin{proposition} \label{1-morphism}
  Let $L=(Y,G,\psi)$ for some cover $Y \rightarrow F$. Assume that patching holds for $G$-torsors.
  Then, the functor $\beta_{L}'''$ is essentially surjective on 1-morphisms.
\end{proposition} 
\begin{proof}
  Let $\G, \G' \in \Gerbe(F,L)$ and let $(\{\G_i \},\{\sigma_k\})$ and $(\{\G'_i \},\{\sigma'_k\})$ denote the images in $\GPAP(\F,L)$.
Given a morphism $(\{\alpha_i\},\{f_k\}) \colon  (\{\G_i \},\{\sigma_k\}) \rightarrow (\{\G'_i \},\{\sigma'_k\})$, we want to construct
a morphism $\G \rightarrow \G'$ whose image in $\GPAP(\F,L)$ is isomorphic to $(\alpha,f)$. \par
Let  $Z \rightarrow F$  be a cover such that $\G(Z) \neq \emptyset \neq \G'(Z)$. 
Thus, $\G|_Z, \G'|_Z \simeq BG$, and therefore $(\{\G_i \},\{\sigma_k\})|_Z=(\{BG_i \},\{\sigma_k|_{Z_k}\})$ and
we can identify $\sigma_k|_{Z_k}$ with the trivial bitorsor $G_k$. The same conclusion holds
for $\G'$ and $(\{\G'_i \},\{\sigma'_k\})$.
Therefore,  over $Z_i$, $\alpha_i$ corresponds to a $G$-bitorsor $P_i$ by Theorem \ref{torsorequi}.
Over $Z_k$, the natural transformation $f_k$  corresponds to the diagram
\begin{center}
  \begin{tikzcd}
    BG_i|_{Z_{k}}  \arrow{r}[name=U]{P_i} \arrow{d}[left]{G_k} & BG_{i}|_{Z_k} \arrow{d}{G_k} \\
    BG_j|_{Z_k}    \arrow{r}[below,name=L]{P_j}  & BG_{j}|_{F_k}  \arrow[Rightarrow,to path=(U) -- (L)\tikztonodes,shorten <= 0.3em,shorten >= 0.3em]{u}{f_{k}}
  \end{tikzcd}
\end{center}
and thus corresponds to an isomorphism of $G_k$-bitorsors $P_i|_{Z_k} \rightarrow P_j|_{Z_K}$.
 By Theorem \ref{bitorsor-patching}, we get a $G$-bitorsor
$P$ defined over $Z$ together with isomorphisms $\phi_i \colon P|_{Z_i}\rightarrow P_i$. This bitorsor in turn defines a morphism
$\alpha \colon \G|_Z \rightarrow \G'|_{Z}$. We claim that it actually descends to a morphism $\G \rightarrow \G'$. According to Lemma \ref{stacks}, we need an
isomorphism of functors $\psi\colon \pr_1^* \alpha \rightarrow \pr_2^* \alpha$ such that 
\begin{center}
  \begin{tikzcd}
\pr_{12}^*\pr_1^* \alpha  \arrow{r}{\pr_{12}^*\psi} \arrow[equal]{d}  & \pr_{12}^*\pr_2^* \alpha  \arrow[equal]{r} & \pr_{23}^*\pr_1^* \alpha \arrow{d}{\pr_{23}^*\psi} \\
\pr_{13}^*\pr_1^* \alpha \arrow{r}{\pr_{13}^*\psi}   & \pr_{13}^*\pr_2^* \alpha \arrow[equal]{r}  & \pr_{23}^*\pr_2^* \alpha
  \end{tikzcd}
\end{center}
commutes. In terms of bitorsors, this means that we need an isomorphism of bitorsors $\psi\colon\pr_1^* P \rightarrow \pr_2^* P$ making the analogous
diagram commute. Such a morphism clearly exists for each $P_i$ as these bitorsors come from morphisms defined over $F_i$.
Furthermore, as these morphisms are compatible with the gluing data, these isomorphisms glue to give a global $\psi\colon\pr_1^* P \rightarrow \pr_2^* P$ by 
Theorem \ref{bitorsor-patching}. Hence, we get a morphism of gerbes $\G \rightarrow \G'$ and it is easy to see that its image
is isomorphic to $(\{\alpha_i\},\{f_k\})$.
\end{proof}

\begin{proposition} \label{2-morphism}
  Let $L=(Y,G,\psi)$ for some cover $Y \rightarrow F$. Assume that patching holds for $G$-torsors.
  Then, the functor $\beta_{L}'''$ is fully faithful on 2-morphisms.
\end{proposition} 
\begin{proof}
  Fix two 1-morphisms of gerbes $\alpha, \beta \colon \G \rightarrow \G'$ and let
$$(\{\alpha_i\},\{f_k\}), (\{\beta_i\},\{g_k\}) \colon (\{\G_i \},\{\sigma_k\}) \rightarrow (\{\G'_i \},\{\sigma'_k\})$$
denote their images in $\GPAP(\F,L)$. \par
Given two 2-morphisms $\psi,\psi' \colon \alpha \rightarrow \beta$ whose image in $\GPAP(\F,L)$ are the same, we want
to prove that $\psi=\psi'$. It is enough to show this after base change, i.e. to prove $\psi|_Z=\psi'|_{Z}$. Thus, we may assume
$\G=\G'=BG$. Thus, $\alpha$ and $\beta$ correspond to $G$-bitorsors $P,Q$ and $\psi,\psi'$ are bitorsor isomorphisms $P \rightarrow Q$.
We then obtain $\psi=\psi'$ immediately from Theorem \ref{bitorsor-patching}. \par
Finally, we need to check fullness. Given a 2-morphism 
$$(\{ \psi_i\}) \colon (\{\alpha_i\},\{f_k\}) \rightarrow (\{\beta_i\},\{g_k\}) $$
we first base change to $Z$. Then, $\psi_i$ correspond to bitorsor isomorphisms and the compatibility condition for
2-morphisms ensures that these isomorphisms glue. Hence, we get $\psi' \colon \alpha|_{Z} \rightarrow \beta|_Z$ by Theorem \ref{bitorsor-patching}.
It remains to show that this morphism descends to a morphism $\psi \colon \alpha \rightarrow \beta$. This follows
from  arguments analogous to the argument to lift the 1-morphism in the proof of Proposition \ref{1-morphism}.
\end{proof}

\begin{theorem} \label{gerbe-patching}
  Let $L=(Y,G,\psi)$ for some cover $Y \rightarrow F$. 
  Assume that  any $L$-gerbe patching problem $(\{\G_i\},\{\nu_k \})$ has property D.
      If patching holds for $G$-torsors, then it also holds for $L$-gerbes, i.e.
  if the functor $\beta_G'$ is an equivalence of 1-categories, then $\beta_{L}'''$ is an equivalence of 2-categories.
\end{theorem}
\begin{proof}
By Proposition \ref{object} and assumption, $\beta_L'''$ is essentially surjective on objects. By Proposition \ref{1-morphism}
$\beta_L'''$ is essentially surjective on 1-morphisms and by Proposition \ref{2-morphism}, $\beta_L'''$  is fully faithful
on 2-morphisms. Hence, $\beta_L'''$ is an equivalence.
\end{proof}

\begin{remark}\label{uniqueness_of_lift}
  Let $L$ be a band and assume that $L$-gerbe patching holds over $\F$. . Let $\G,\G'$ be $L$-gerbes over $\F$ such that
  $\beta'''_{L,\F}(\G)$ and $\beta'''_{L,\F}(\G')$ are isomorphic.
  Then, there is an equivalence $\G \simeq \G'$, unique up to unique isomorphism, as $\beta'''_{L,\F}$ is an equivalence
  of 2-categories.
\end{remark}
We will now conclude this section with a sufficient condition for the technical assumption of the
above theorem to be satisfied.

\begin{proposition}\label{technical}
  Let $L=(Y,G,\psi)$ be a band over $F$. Assume that patching holds for $\HH^2(\bullet,L)$ and
  that bitorsor factorization holds for any cover $Z \rightarrow Y$. Then,  any $L$-gerbe patching
  has property D.
\end{proposition} 
\begin{proof}
  Let $\alpha_i \in \HH^2(F_i,L_i)$ denote the equivalence class of $\G_i$. Then, we have
  that $\alpha_i=\alpha_j \in \HH^2(F_k,L_k)$ by assumption. Hence, there is some
  $\alpha \in \HH^2(F,L)$ inducing $\alpha_i$ for all $i \in I_v$. Let $\G$ be an $L$-gerbe
  representing $\alpha$ and let $Z \rightarrow Y$ be a cover such that $\G(Z) \neq \emptyset$.
  Then, $\G_i(Z_i) \neq \emptyset$ by construction. Fix equivalences $\G_i \rightarrow BG_i$.
  Then, for $(i,j,k) \in S_I$, we have a morphism 
  \begin{center}
    \begin{tikzcd}[column sep=small, row sep=small]
BG_i  \arrow{d}    & BG_j \\
\G_i    \arrow{r}{\sigma_k}  & \G_j \arrow{u}
    \end{tikzcd}
  \end{center}
defined over $Z_k$.
Let $P$ be the $G_k$-bitorsor corresponding to the composition $BG_i \rightarrow BG_j$.
By assumption, we  there is a $G_i$-bitorsor $P_i$ over $Z_i$ and a $G_j$-bitorsor $G_j$ over $Z_j$
such that $P\simeq P_i \wedge^{G_k} P_j$ over $Z_k$. Let $x_i \in \G_i(Z_i)$ correspond
to $P^{\op}_i$ with respect to the chosen equivalence $\G_i \rightarrow BG_i$ and let $x_j$ correspond
to $P_j$. Then, it is straightforward to check that $\sigma_k(x_i|_{Z_k})$ and $x_j|_{Z_k}$
are isomorphic in $\G_j(Z_k)$.
\end{proof}

\begin{corollary}
  Under the assumptions of Proposition \ref{technical},
  $L$-gerbe patching holds over $\F$.
\end{corollary}

\subsection{A Mayer-Vietoris sequence and a local-global principle for gerbes} \label{gerbe-lg}

In this section, we will fix a finite inverse factorization system $\F$ with inverse limit $F$.
We will also fix a  group sheaf $G$ defined over $F$ for which $G$-torsor patching holds.
 Let $L$ be the band induced by $G$.

\begin{defn}
 We say that $G$-gerbes satisfy the local-global principle
with respect to $\F$ if  for any $G$-gerbe $\G$ over $F$,
  \begin{align*}
    \G_i \simeq BG|_{F_i}
  \end{align*}
  for all $i \in I_v$ implies 
  \begin{align*}
    \G \simeq BG.
  \end{align*}
\end{defn}

Recall that $G$-bitorsors over $F$ are classified by $\HH^0(F,G \rightarrow \Aut(G))$
and $G$ gerbes are classified by $\HH^1(F, G \rightarrow \Aut(G))$ (cf. Appendix \ref{hyper}).

Note that we obtain two maps
\begin{align*}
  \prod_{i \in I_v} \HH^1(F_i, G \rightarrow \Aut(G)) \rightrightarrows
    \prod_{k \in I_e} \HH^1(F_k, G \rightarrow \Aut(G)).
\end{align*}
via base change.

\begin{lemma}
  There is a map of pointed sets
  \begin{align*}
    \prod_{k \in I_k} \HH^0(F_k, G \rightarrow \Aut(G)) \rightarrow
      \HH^1(F, G \rightarrow \Aut(G))
  \end{align*}
\end{lemma}
\begin{proof}
   We can define the map 
as follows: Given $G$-bitorsors $P_k $ over $F_k$, consider the $G$-gerbe patching
problem $(\{ BG|_{F_i}\}_i, \{ P_k \}_k )$. By Theorem \ref{gerbe-patching}, there is a $G$-gerbe
$\G$ over $F$ in the essential preimage of the patching problem.
We map the equivalence class of $(P_k)_k$ to the equivalence class $\G$. Note that this
is independent of the choice of $\G$ (cf. Remark \ref{uniqueness_of_lift}).
\end{proof}
 
We can put these maps in a Mayer-Vietoris type sequence.
\begin{theorem}[Mayer-Vietoris for non-abelian hypercohomology (2)] \label{mayer2}
  Assume patching holds for $G$-torsors and $G$-gerbes over $\F$.
  Then, there is an exact sequence 
   \begin{adjustwidth}{-2em}{-2em}
  \begin{center}
    \begin{tikzcd}[row sep =small, column sep = small]
   1 \arrow{r} &    \HH^{-1}(F, G \rightarrow \Aut(G)) \arrow{r}&
    \prod_{i \in I_v} \HH^{-1}(F_i, G \rightarrow \Aut(G))\arrow{r} \arrow[draw=none]{d}[name=Z, shape=coordinate]{} & 
    \prod_{k \in I_e} \HH^{-1}(F_k, G \rightarrow \Aut(G))  \arrow[rounded corners,to path={ -- ([xshift=2ex]\tikztostart.east)|- (Z) [near end]\tikztonodes-| ([xshift=-2ex]\tikztotarget.west)-- (\tikztotarget)}]{dll} \\
 \;&     \HH^0(F, G \rightarrow \Aut(G)) \arrow{r}&
    \prod_{i \in I_v} \HH^0(F_i, G \rightarrow \Aut(G))\arrow{r} \arrow[draw=none]{d}[name=Z, shape=coordinate]{} & 
    \prod_{k \in I_e} \HH^0(F_k, G \rightarrow \Aut(G))  \arrow[rounded corners,to path={ -- ([xshift=2ex]\tikztostart.east)|- (Z) [near end]\tikztonodes-| ([xshift=-2ex]\tikztotarget.west)-- (\tikztotarget)}]{dll} \\
  \; &    \HH^1(F, G \rightarrow \Aut(G)) \arrow{r}& 
    \prod_{i \in I_v} \HH^1(F_i, G \rightarrow \Aut(G))\arrow[yshift=3]{r} \arrow[yshift=-3]{r}&
    \prod_{k \in I_e} \HH^1(F_k, G \rightarrow \Aut(G))
    \end{tikzcd}
  \end{center}
  \end{adjustwidth}
\end{theorem} 
\begin{proof}
  The exactness in the first two rows is the content of  Theorem \ref{mayer1}.
  The exactness at $ \prod_{i \in I_v} \HH^1(F_i, G \rightarrow \Aut(G))$ follows
  from  gerbe  patching (Theorem \ref{gerbe-patching}).
  The exactness at  $\HH^1(F, G \rightarrow \Aut(G))$ follows immediately from Theorem \ref{morita}.
\end{proof}

From this exact sequence, we can deduce a necessary and sufficient criterion for the local-global
principle for gerbes to hold in terms of bitorsor factorization.

\begin{theorem} \label{local-global1}
Assume that gerbe patching holds for $L$-gerbes. Then, $L$-gerbes satisfy a local-global
principle with respect to patches if and only $G$ satisfies bitorsor factorization.  
\end{theorem}

\begin{remark}
  These results are analogous to the results in \cite{Harbater2011} concerning
  local-global principles for $G$-torsors. They prove that local-global principle
  for $\HH^1(F,G)$ is equivalent to factorization of $\prod_k \HH^0(F_k,G)$. In other words, local-global
  principle for $G$-torsors is equivalent to $G$ satisfying factorization.
\end{remark}

\section{Patching over arithmetic curves} \label{curve}

We will now apply the results on bitorsor and gerbe patching in the case of
patching
over arithmetic curves, (cf. Section \ref{arithmetic}).
Let $T$ be a complete discretely valued ring with
field of fraction $K$, uniformizer $t$ and residue field $k$.
Let $\Xhat$ be a projective, integral and normal $T$-curve with function field $F$ and let $X$ denote its closed fiber.
For a non-empty set of closed  points $\calP\subset X$  including all points where distinct irreducible components meet,
let $\F$ denote a finite inverse factorization system with inverse limit $F$ as introduced in Section \ref{arithmetic}.

\subsection{Patching and local-global principle for bitorsors}
Let $G,H$ be a linear algebraic groups over $F$. 

Recall the natural functor
\begin{align*}
  \beta_{(G,H)}'' \colon \Bitorsor(G,H)(F) \rightarrow  \BPAP(G,H)(\F).
\end{align*}
introduced in Section \ref{sec:bit-pat}.

\begin{theorem} \label{curve-bitorsor}
  Let $G,H$ be linear algebraic groups over $F$. Then,
  patching holds for $(G,H)$-bitorsors over $\F$, i.e.
  the functor $\beta_{(G,H)}''$ is an equivalence
\end{theorem} 
\begin{proof}
  Follows from Theorem \ref{curve-torsor} and Theorem \ref{bitorsor-patching}.
\end{proof}

As a corollary, we immediately obtain a criterion for when $G$-bitorsors
satisfy a local-global principle.

\begin{corollary} \label{curve-local-bit}
  The local-global principle for $G$-bitorsors holds over $\F$ iff
  $Z(G)$ satisfies factorization over $\F$.
\end{corollary} 
\begin{proof}
  Follows from Corollary \ref{local-global-bitorsor} and Theorem \ref{curve-bitorsor}.
\end{proof}

\begin{theorem}
  Let $G$ be a linear algebraic group whose center $Z$ is rational. Then,
  $G$-bitorsors satisfy local-global principle over $\F$ iff $Z$ is connected or
  $\Gamma$ is a tree.
\end{theorem} 
\begin{proof}
  In \cite[Corollary 6.5]{Harbater2011}, Harbater, Hartmann and Krashen proved that a rational linear algebraic group
  $H$ satisfies factorization over $\F$ iff $H$ is connected or $\Gamma$ is a tree.
  Hence, the result follows from  Corollary \ref{curve-local-bit}.
\end{proof}

\subsection{A Mayer-Vietoris sequence and gerbe patching over curves}

We will now investigate, when gerbe patching holds over arithmetic curves.
Let us first collect some results related to the technical assumption in Theorem \ref{gerbe-patching}.

\begin{theorem}[\cite{Colliot-Thelene2017}] \label{suresh}
  Assume that $\chara(K)=0=\chara(k)$. Let $(L_i)_{i \in I}$ be a collection of finite
  separable field extensions $L_i / F_i$. Then, there is a finite separable field extension
  $E/F$ such that $E_i$ dominates $L_i$.
\end{theorem}

\begin{corollary} \label{assum-gerbe-char0}
  Assume that $\chara(k)=0$. Then, every $L$-gerbe patching problem $$(\{\G_{i}\}_{i \in I_v}, \{\nu_k \}_{k \in I_e})$$
   over $\F$ has property D.
\end{corollary} 
\begin{proof}
  Pick covers $Z_i \rightarrow F_i$ for $i \in I_v$ and $z_i \in \G_i(Z_i)$. Note that $Z_i$ is a product
  of finite separable field extensions of $F_i$. By Theorem \ref{suresh}, there is a cover 
  $Z' \rightarrow F$ such that $Z'_i$ dominates $Z$.  \par
  While $\nu_k(x_i|_{Z'_k})$ and $x_j|_{Z'_k}$ may not be isomorphic in $\G_j(Z'_k)$, they are locally isomorphic,
  so there are covers $Y_k \rightarrow Z'_k$ such $\nu_k(x_i|_{Y_k})$ and $x_j|_{Y_k}$ are isomorphic.
  Again using Theorem \ref{suresh}, there is a cover $Y \rightarrow F$ dominating $Z \rightarrow F$
  such that $Y'_k \rightarrow F$ dominates $Y_k \rightarrow F_k$.  \par
  Then, the choice $Y' \rightarrow F$ and $x_i|_{Y'_i}$ proves the claim.
\end{proof}

In the case where $\chara(k)=p$, we want to use Proposition \ref{technical} to prove that
every gerbe patching problem has property D. The next proposition proves that this is true under
some mild assumptions.

\begin{proposition} \label{assum-gerbe-charp}
  Let $L$ be a band over $F$ such that $Z(L)$ is a linear algebraic group over $F$
  with finite order not divided by $\chara(k)$. Then, if $\HH^2(F,L) \neq \emptyset$, 
  patching holds for $\HH^2(\bullet, L)$ over
  $\F$.
\end{proposition} 
\begin{proof}
  By  \cite[Theorem 3.1.3]{Harbater2012a}, patching holds for $\HH^2(\bullet, Z(L))$.
  Thus, the conclusion follows by Proposition \ref{patching-coho}.
\end{proof}

We are now ready to state our main result on gerbe patching over arithmetic curves.

\begin{theorem} \label{curve-gerbes}
  Let $L=(Y,G,\psi)$ be a band over $F$ with $G$ being a linear algebraic group over $Y$.
  Assume either  
  \begin{itemize}
    \item $\chara(k)=0$  or
    \item $\chara(k)=p>0$ and  $Z(L)$ has finite order not divisible by $\chara(k)$ and
      $G$-bitorsor factorization holds over $\F_Z$  for every cover $Z \rightarrow Y$.
  \end{itemize}
  Then, gerbe patching holds over $\F$, i.e. the functor
  $\beta'''_{L}$ is a 2-equivalence.
\end{theorem} 
\begin{proof}
  If $\chara(k)=0$, this follows from Theorem \ref{curve-torsor}, Theorem \ref{gerbe-patching} and Theorem \ref{assum-gerbe-char0}.
  If $\chara(k)=p$, then it follows from Theorem \ref{curve-torsor}, Theorem \ref{gerbe-patching} and Proposition \ref{assum-gerbe-charp}.
\end{proof}

\begin{corollary}[Mayer-Vietoris of non-abelian hypercohomology over curves] \label{mayer-curves}
  Let $G$ be a linear algebraic group defined over $F$ and let $L$ denote the
  associated band.
  Under the assumption of Theorem \ref{curve-gerbes}, there is an exact sequence
  of pointed sets
   \begin{adjustwidth}{-2em}{-2em}
   \begin{center}
    \begin{tikzcd}[row sep =small, column sep = small]
  1 \arrow{r} &     \HH^{-1}(F, G \rightarrow \Aut(G)) \arrow{r}&
    \prod_{i \in I_v} \HH^{-1}(F_i, G \rightarrow \Aut(G))\arrow{r} \arrow[draw=none]{d}[name=Z, shape=coordinate]{} & 
    \prod_{k \in I_e} \HH^{-1}(F_k, G \rightarrow \Aut(G))  \arrow[rounded corners,to path={ -- ([xshift=2ex]\tikztostart.east)|- (Z) [near end]\tikztonodes-| ([xshift=-2ex]\tikztotarget.west)-- (\tikztotarget)}]{dll} \\
   \; &   \HH^0(F, G \rightarrow \Aut(G)) \arrow{r}&
    \prod_{i \in I_v} \HH^0(F_i, G \rightarrow \Aut(G))\arrow{r} \arrow[draw=none]{d}[name=Z, shape=coordinate]{} & 
    \prod_{k \in I_e} \HH^0(F_k, G \rightarrow \Aut(G))  \arrow[rounded corners,to path={ -- ([xshift=2ex]\tikztostart.east)|- (Z) [near end]\tikztonodes-| ([xshift=-2ex]\tikztotarget.west)-- (\tikztotarget)}]{dll} \\
 \; &     \HH^1(F, G \rightarrow \Aut(G)) \arrow{r}& 
    \prod_{i \in I_v} \HH^1(F_i, G \rightarrow \Aut(G))\arrow[yshift=3]{r} \arrow[yshift=-3]{r}&
    \prod_{k \in I_e} \HH^1(F_k, G \rightarrow \Aut(G))
    \end{tikzcd}
  \end{center}
  \end{adjustwidth}
\end{corollary} 
\begin{proof}
  Follows from Theorem \ref{mayer2} and Theorem \ref{curve-gerbes}.
\end{proof}

\subsection{Factorization of bitorsors over curves}

We will now investigate which group schemes $G$ over $F$ admit $G$-bitorsor
factorization over $\F$.  Let $\F$ be indexed by $I=I_v \sqcup I_e$ with associated graph $\Gamma$ (cf. Section \ref{patching} for details).
The short exact sequence $1 \rightarrow (1 \rightarrow \Aut(G)) \rightarrow (G \rightarrow \Aut(G)) \rightarrow
    (G \rightarrow 1) \rightarrow 1$ of crossed modules induces the long exact sequence
\begin{equation} \label{les}
\begin{tikzcd}[column sep=small]
  1 \arrow{r} &\HH^{0}(Z(G)) \arrow{r} & \HH^{0}(G) \arrow{r}  \arrow[draw=none]{d}[name=Z, shape=coordinate]{} 
  &  \HH^{0}(\Aut(G))  \arrow[rounded corners,to path={ -- ([xshift=2ex]\tikztostart.east)|- (Z) [near end]\tikztonodes-| ([xshift=-2ex]\tikztotarget.west)-- (\tikztotarget)}]{dlll} \\
    \HH^{0}(G \rightarrow \Aut(G))\arrow{r} & \HH^{1}(G) \arrow{r} & \HH^{1}(\Aut(G))    \arrow{r}  & \HH^{1}(G \rightarrow \Aut(G)) & \; 
\end{tikzcd}
\end{equation}
(compare Example \ref{myles}).

We will first consider the case of finite constant group schemes with trivial center. 
\begin{theorem}\label{finite}
  Assume that $G$ is a finite constant group scheme over $F$ with trivial center
  and that $\Gamma$ is a tree.
  Then, $G$ satisfies bitorsor factorization over $\F$.
\end{theorem} 
\begin{proof}
  Given a collection $\{P_b \}_{b \in \calB}$ with $P_b$ a $G$-bitorsor over $F_b$, we need to show that
  there are $\{P_U \}_{U \in \calU}$ and $\{P_p \}_{p \in \calP}$ such that $P_U$ is a $G$-bitorsor
  over $F_U$, $P_p$ is a $G$-bitorsors over $F_P$ and $P_U|_{F_b} \wedge^G P_p|_{F_b} \simeq F_b$
  whenever $b$ is a branch at $U$ and $p$.

  By Theorem \ref{hyperbit}, this is equivalent to showing that
  \begin{align*}
  \prod_{U \in \calU}\HH^{0}(F_U,G \rightarrow \Aut(G))\times  \prod_{p \in \calP}\HH^{0}(F_p,G \rightarrow \Aut(G)) \rightarrow \prod_{b \in \calB}\HH^{0}(F_{b}, G \rightarrow \Aut(G))
  \end{align*}
  is surjective.
  Since $Z(G)=\{e\}$, the sequence
  \begin{align*}
    1 \rightarrow G \rightarrow \Aut(G) \rightarrow \Aut(G) /G \rightarrow 1
  \end{align*}
  is exact. Since $G$ is constant, so is $\Aut(G)$ and the long exact sequence associated
  to the short exact sequence above reads
  \begin{align*}
    1 \rightarrow \HH^1(F_{b},G) \rightarrow \HH^1(F_{b},\Aut(G)) \rightarrow \ldots
  \end{align*}
  for any $b \in \calB$. Since $\HH^1(F_{b},G) \rightarrow \HH^1(F_{b},\Aut(G)$ is injective, it follows
  from sequence (\ref{les})
  that $\HH^0(F_b,\Aut(G) \rightarrow \HH^0(F_b,G \rightarrow \Aut(G))$ is surjective.
  It is thus enough to show that
  \begin{align*}
    \prod_{U \in \calU}\HH^{0}(F_U, \Aut(G))\times  \prod_{p \in \calP}\HH^{0}(F_p, \Aut(G)) \rightarrow \prod_{b \in \calB} \HH^{0}(F_{b}, \Aut(G))
  \end{align*}
  is surjective,i.e. that $\Aut(G)$ satisfies factorization. Since $\Aut(G)$ is also a finite, constant
  group scheme, it satisfies factorization by assumption and Theorem \ref{simfact}.
\end{proof}

\begin{examples}
Examples of group schemes satisfying the assumptions of theorem \ref{finite}
include  $S_n$ for any $n>2$, $A_n$ for $n>3$ and any finite, nonabelian simple group.
\end{examples}

\begin{theorem} \label{factor2}
  Let $G$ be an algebraic group over $F$ such that the natural map $G \rightarrow \Aut(G)$
  is an isomorphism. If $G$ satisfies factorization over $\F$, then $G$ satisfies bitorsor factorization
  over $\F$.
\end{theorem} 
\begin{proof}
  Using the assumption, we can see that the map $\HH^0(F,G)=\HH^0(F,\Aut(G)) \rightarrow \HH^0(F, G \rightarrow \Aut(G))$
  is surjective. Hence, the claim follows immediately.
\end{proof}

\begin{examples} \label{noouter}
Let $G$ be semisimple adjoint, Then, $G \rightarrow \Aut(G)$ is an isomorphism if  $G$ is of type $A_1$, $B_n$, $C_n$, 
$E_7$, $E_8$, $F_4$, and $G_2$ (compare 24.A and \cite[Proposition 25.15]{Knus}). \par
Thus, if $G$ satisfies factorization over $\F$, then it also satisfies bitorsor factorization.
\end{examples}

Let now $G$ be a semisimple group whose adjoint group admits no outer automorphism.
Then, we have a short exact sequence of crossed modules (see Appendix \ref{hyper} for the definition and see Corollary 25.17 in \cite{Knus}
for exactness):
\begin{align*}
  1 \rightarrow (Z \rightarrow 1) \rightarrow (G \rightarrow \Aut(G)) \rightarrow
    (G/Z \rightarrow \Aut(G/Z)) \rightarrow 1
\end{align*}
We thus obtain the following exact sequence of hypercohomology groups.
\begin{equation} \label{les2}
\begin{tikzcd}[column sep=small]
  1\arrow{r}&  \HH^0(F,Z) \arrow{r}  &  \HH^0(F,Z) \arrow{r} \arrow[draw=none]{d}[name=Z, shape=coordinate]{}  &  1 \arrow[rounded corners,to path={ -- ([xshift=2ex]\tikztostart.east)|- (Z) [near end]\tikztonodes-| ([xshift=-2ex]\tikztotarget.west)-- (\tikztotarget)}]{dll} \\
 & \HH^1(F,Z)\arrow{r}& \HH^0(F,G \rightarrow \Aut(G)) \arrow{r} \arrow[draw=none]{d}[name=Z2, shape=coordinate]{}&  \HH^0 (F,G/Z \rightarrow \Aut(G/Z)) \arrow[rounded corners,to path={ -- ([xshift=2ex]\tikztostart.east)|- (Z2) [near end]\tikztonodes-| ([xshift=-2ex]\tikztotarget.west)-- (\tikztotarget)}]{dll} \\
  & \HH^2(F,Z) \arrow{r} &\HH^1(F,G \rightarrow \Aut(G))& 
\end{tikzcd}
\end{equation}

\begin{lemma} \label{0map}
  Let $G$ be a semisimple group whose adjoint group $G/Z$ admits no outer automorphisms.
  Then, the map
  \begin{align*}
    \HH^0(F,G/Z \rightarrow \Aut(G/Z)) \rightarrow \HH^2(F,Z)
  \end{align*}
  from the exact sequence (\ref{les2}) is the zero map.
\end{lemma}
\begin{proof}
  Consider the diagram 
  \begin{center}
    \begin{tikzcd}[column sep=small, row sep=small]
      \;      & \HH^1(F,G/Z) \arrow{rd}{\delta} & \; \\
      \HH^0(F,G/Z \rightarrow \Aut(G/Z))  \arrow{rr} \arrow{ru}    & \; & \HH^2(F,Z).
    \end{tikzcd}
  \end{center}
  Note that this diagram does not commute: a cocycle $(f,\lambda) \in Z^0(F,G \rightarrow \Aut(G))$
  maps to $\lambda( \delta(f))$ as opposed to $\delta(f)$. However, as $\delta(f)=0$ implies
  $\lambda(\delta(f))=0$, it is enough to show that  
  \begin{align*}
    \HH^0(F,G/Z \rightarrow \Aut(G/Z)) \rightarrow \HH^1(F,G/Z) 
  \end{align*}
  is the zero map. \par
  From the exact sequence (\ref{les}), we obtain that 
  \begin{align*}
       \HH^0(F,G/Z \rightarrow \Aut(G/Z)) \rightarrow \HH^1(F,G/Z) \rightarrow \HH^1(F,\Aut(G/Z))                                             
  \end{align*}
  is exact. But, as $G/Z$ is isomorphic to $\Aut(G/Z)$ by assumption, the claim follows.
\end{proof}

\begin{theorem} \label{center}
  Let $G$ be a semisimple group such that $G/Z$ admits no outer automorphisms. If $G/Z$ and $Z$ satisfy bitorsor factorization over $\F$,
  then so does $G$.
\end{theorem}
\begin{proof}
Let us start with $\alpha \in \HH^0(F_{b},G \rightarrow \Aut(G))$. Let
$\beta$ denote its image in $\HH^1(F_{b},G/Z \rightarrow \Aut(G/Z))$.
By assumption, there are $\beta_U \in \HH^0(F_U,G/Z \rightarrow \Aut(G/Z))$
and $\beta_p \in \HH^0(F_p, G/Z \rightarrow \Aut(G/Z))$ such that
$\beta_U \beta_p=\beta_{b} \in \HH^0(F_{b},G/Z \rightarrow \Aut(G/Z))$.
By Sequence (\ref{les2}) and  Lemma \ref{0map}, we can lift $\beta_U$ and $\beta_P$ to elements $\tilde{\alpha}_U$ and $\tilde{\alpha}_p$ in
$\HH^0(F_U,G \rightarrow \Aut(G))$ and $\HH^0(F_p, G \rightarrow \Aut(G))$ respectively. Consider the element
$\tilde{\alpha}=\tilde{\alpha}^{-1}_U \alpha \tilde{\alpha}^{-1}_p \in \HH^0(F_{b}, G \rightarrow \Aut(G))$. By construction,
$\tilde{\alpha}$ comes from an element $\gamma \in \HH^1(F_{b},Z)$. By assumption, there are
$\gamma_U \in \HH^1(F_U,Z)$ and $\gamma_p \in \HH^1(F_p,Z)$ satisfying $\gamma_U\gamma_p=\gamma \in \HH^1(F_{b},Z)$.
Consider the elements $\alpha_U=\gamma_U \tilde{\alpha}_u \in \HH^0(F_U, G \rightarrow \Aut(G))$ and
 $\alpha_p=\gamma_p \tilde{\alpha}_p \in \HH^0(F_p, G \rightarrow \Aut(G))$. Then, $\alpha=\alpha_U \alpha_p$. This follows from
\begin{align*}
  \alpha_U^{-1} \alpha \alpha_p^{-1} &= \gamma_U^{-1} \tilde{\alpha}_U^{-1} \alpha \tilde{\alpha}_p^{-1} \gamma_p^{-1} \\
                                   &= \gamma_U^{-1} \gamma \gamma_p^{-1} \\
                                   &= 1.
\end{align*}
\end{proof}

We will now turn our attention to groups of type $A_n$. For this, let $D$ be a central simple algebra
over $F$.
\begin{lemma}\label{autsl1}
  The group $\Aut(\SL_1(D))$ satisfies factorization over $\F$ if $\Gamma$ is a tree.
\end{lemma}
\begin{proof}
  Note that we have a short exact sequence
  \begin{align*}
    1 \rightarrow \PGL_1(D)(K) \rightarrow \Aut(\SL_1(D))(K) \rightarrow \Z / 2\Z \rightarrow 1
  \end{align*}
  for any field $K/F$. The result now follows from noting that $\PGL_1(D)$ and $\Z/2\Z$ satisfy factorization
  over $\F$ as they are both rational and $\Gamma$ is a tree (see Theorem \ref{simfact}).
\end{proof}

\begin{proposition}\label{trivial}
  Let $D$ be a central simple algebra over $K$. Then, the map of pointed sets
  \begin{align*}
    \HH^1(K,\SL_1(D) \rightarrow \HH^1(K,\Aut(\SL_1(D)))
  \end{align*}
 has trivial image. In particular,  the map
  $\HH^0(\SL_1(D) \rightarrow \Aut(\SL_1(D))) \rightarrow \HH^1(\SL_1(D))$ is surjective.
\end{proposition} 
\begin{proof}
Note that a piece of sequence (\ref{les}) for $G= \SL_1(D)$ reads
\begin{align} \label{sldles}
  \HH^0(\Aut(\SL_1(D)) \rightarrow \HH^0(\SL_1(D) \rightarrow \Aut(\SL_1(D))) \rightarrow \HH^1(\SL_1(D)) \rightarrow \HH^1(\Aut(\SL_1(D))).
\end{align}
Thus, injectivity of $\HH^1(K,\SL_1(D) \rightarrow \HH^1(K,\Aut(\SL_1(D)))$ implies surjectivity of
 $\HH^0(\SL_1(D) \rightarrow \Aut(\SL_1(D))) \rightarrow \HH^1(\SL_1(D))$.
  Note that the map $\SL_1(D) \rightarrow \Aut(\SL_1(D))$ factors through
  $\SL_1(D) \rightarrow \GL_1(D)$. Hence, $\HH^1(K,\SL_1(D) \rightarrow \HH^1(K,\Aut(\SL_1(D)))$ factors through
  $\HH^1(K,\SL_1(D) \rightarrow \HH^1(K,\GL_1(D))$.
  The claim thus follows from $\HH^1(K, \GL_1(D))=0$.
\end{proof}

\begin{theorem}\label{sld}
  Let $D$ be a a central simple algebra over $F$ and let $\Gamma$ be a tree. Then,
  $\SL_1(D)$ satisfies bitorsor factorization.
\end{theorem} 
\begin{proof}
  Fix a collection of
  $$\alpha_{b} \in \HH^0(F_{b},\SL_1(D) \rightarrow \Aut(\SL_1(D)))$$
  for $b \in \calB$
  and let $\beta_{b}$ denote their images
  in $\HH^1(F_{b},\SL_1(D))$ (along sequence \ref{sldles}) when $b$ corresponds to a branch at $p$. Recall that
  $$\HH^1(\F_{b},\SL_1(D))= F^*_{b} / \Nrd(D_{b})$$
  holds. Let $n$ be the index of $D$. Since $F_b$ is a completion of $F_p$ at a discrete valuation,
  the map $F_p^* \rightarrow F_b^* / F_b^{*n}$ is surjective. As $F^{*n}_n \subset \Nrd(D^*_b)$,
  it follows that $F_p^* \rightarrow F_b^* / \Nrd(D^*_b)$ is surjective whenever $b$ is a branch
  at $p$. Thus, by weak approximation, there are $\beta_p \in \HH^1(F_p, \SL_1(D))$ for all $p \in \calP$ such that
  $\beta_p=\beta_b \in \HH^1(\F_{b},\SL_1(D))$ whenever $b$ is a branch at $p$.
  
  By Proposition \ref{trivial}, there are $\alpha_p \in \HH^0(F_p, \SL_1(D) \rightarrow \Aut(\SL_1(D)))$ mapping onto $\beta_p$ for all $p \in \calP$.
  Let now $b$ be a branch at $p$ and $U$. Then,
  $\alpha_b \alpha_p^{-1}$ maps onto $0$ in $\HH^1(F_{b},\SL_1(D))$ by construction. Hence, there is $\nu_{b} \in \HH^0(F_{b},\Aut(\SL_1(D_{b})))$ mapping onto 
  $\alpha_b \alpha_p^{-1}$.
  By Lemma \ref{autsl1}, there exist 
  $\nu_U  \in \HH^0(F_{U},\Aut(\SL_1(D)))$ and $ \nu_p \in \HH^0(F_{p},\Aut(\SL_1(D)))$
  such that their product in $\HH^0(F_{b},\Aut(\SL_1(D)))$ is $\nu_{b}$ for all $(U,p,b) \in S_I$. \par
  Let $\tau_U$ and $\tau_p$ denote the images of $\nu_U$ and $\nu_p$ in $\HH^0(F_{U},\SL_1(D) \rightarrow \Aut(\SL_1(D))))$ 
  and  $\HH^0(F_{p},\SL_1(D) \rightarrow \Aut(\SL_1(D)))$ respectively.
  Then, the collection of $\{ \tau_U\}_{U \in \calU}$ and $\{\tau_p \alpha_p\}_{p \in \calP}$  give a factorization of $\{\alpha_b\}_{b \in \calB}$.
\end{proof}

\subsection{Local-global principles for gerbes}

Building on our results on bitorsor factorization, we now use Theorem \ref{local-global1}
to obtain local-global principles for gerbes.

\begin{theorem} \label{curve_local}
  Let $G$ be a linear algebraic group over $F$ with center $Z$. Then, the local global
  principle for $G$-gerbes with respect to patching holds if 
  \begin{itemize}
    \item  $\chara(k)=0$ and one of the following hold:
      \begin{itemize}
        \item $\Gamma$ is a tree and $G$ is a finite constant group scheme with trivial center,
        \item $G$ is connected, rational, semisimple, adjoint of type $A_1,B_n,C_n,E_7,E_8,F_4$ or $G_4$,
        \item $G$ is semisimple such that $G/Z$ admits no outer automorphism and $Z,G/Z$ satisfy bitorsor factorization,
        \item $\Gamma$ is a tree and $G=\SL_1(D)$ where $D$ is a central simple algebra over $F$,
      \end{itemize}
    \item  $\chara(k)=p>0$, $Z$ has finite order not divided by $\chara(k)$ and $G$ and $\Gamma$ are as in the case of $\chara(k)=0$.
  \end{itemize}
\end{theorem} 
\begin{proof}
  All results use Theorem \ref{gerbe-patching} and Theorem \ref{local-global1}.
  The results follow (in order) from  Theorem \ref{finite}, Example \ref{noouter} and  \cite[Theorem 3.6]{Harbater2009}, Theorem \ref{center},
  and Theorem \ref{sld}.
\end{proof}

\subsection{Local-global principle for homogeneous spaces}
In this section, we will apply our local-global principle for gerbes to deduce
local-global principles for points on homogeneous spaces.

Let $H$ be a linear algebraic group and let $X$ be a homogeneous space under $H$ over $F$.
Assume that the stabilizer  of a geometric point $x \in X(F^{\sep})$ is isomorphic
to $G|_{F^{\sep}}$ for some linear algebraic group $G \subset H$ defined over $F$.
Consider the quotient stack $[X/H]$. Objects of this stack are given by tuples
$(Y,P,\phi)$ where $Y$ is an $F$-scheme, $P$ is a principal homogeneous space under 
$H|_{Y}$ over $Y$ and $\phi\colon P \rightarrow X|_Y$ is an $H$-equivariant map.
The stack $[X/H]$ is a $G$-gerbe and it is trivial if and only if there is a $H$-equivariant map
$P \rightarrow X$ where $P$ is principal homogeneous space under $H$ over $F$ (cf. \cite[Remark 7.7.1]{Borovoi1993}).

Let $Y$ be an $F$-scheme. If $X(Y)\neq \emptyset$, pick a point $x \in X(Y)$.
This defines an $H$-equivariant morphism  $\phi_x \colon H|_Y \rightarrow X|_Y$
 which gives an element $(Y,H|_Y,\phi_x)\in [X/H](Y)$. 
Hence, $X(Y) \neq \emptyset$ implies$[X/H](Y) \neq \emptyset$. 
Assume now that  $H$ is special (i.e. $\HH^1(K,H)=\{e\}$ for all field extensions $K/F$).
Let $K/F$ be a field extension. If $[X/H](K)\neq \emptyset$, then there is a principal homogeneous space $P$
under $H|_K$ over $K$ together with an $H$ equivariant map $\phi \colon P \rightarrow X$. As $H$ is special,
$P$ admits a point $p \in P(K)$. Hence, $\phi(p)\in X(K)$. Thus, $[X/H](K) \neq \emptyset$ implies
$X(K) \neq \emptyset$.

\begin{defn}
  We say that the local-global principle holds for $X$ over $\F$ if
  $X(F_i) \neq \emptyset$ for all $i \in I_v$ implies $X(F) \neq \emptyset$.
\end{defn}

\begin{theorem} \label{homspace}
  Let $H$ be special and let $X$ be a homogeneous space under $H$. Let $G \subset H$ be a linear algebraic group
  over $F$ such that $G|_{F^{\sep}}$  is isomorphic to the stabilizer of a geometric
  point of $X$. If $\chara(k)=0$, the local-global principle for $X$ holds if and only if
  bitorsor factorization holds for $G$. \par
  If $\chara(k)=p>0$,  assume that 
  \begin{itemize}
    \item $p \nmid |Z(G)| < \infty$,
    \item $G$-bitorsor factorization holds over $\F_Z$ for every cover $Z \rightarrow F$.
  \end{itemize} 
  Then, the local-global principle holds for $X$ over $\F$.
\end{theorem} 
\begin{proof}
  Since $H$ is special, we have $X(F) \neq \emptyset$ iff $[X/H](F) \neq \emptyset$ as well as
  $X(F_i) \neq \emptyset$ iff $[X/H](F_i) \neq \emptyset$ by the discussion above. Hence, the local-global principle
  for $X$ holds iff it holds for $[X/H]$. The result follows from Theorem \ref{local-global1}
  and Theorem \ref{curve-gerbes}.
\end{proof}

\begin{corollary} \label{above}
  Let $H$ be special (e.g. $H=\SL_n$ or $H=\SP_{2n}$)  and let $X$ be a homogeneous space under $H$. Let $G \subset H$ be a linear algebraic group
  over $F$ such that $G|_{F^{\sep}}$  is isomorphic to the stabilizer of a geometric
  point of $X$.
   \begin{itemize}
    \item If $\chara(k)=0$,  assume that  one of the following holds: 
      \begin{itemize}
        \item $\Gamma$ is a tree and $G$ is a finite constant group scheme with trivial center,
        \item $G$ is connected, rational, semisimple, adjoint of type $A_1,B_n,C_n,E_7,E_8,F_4$ or $G_4$,
        \item $G$ is semisimple such that $G/Z$ admits no outer automorphism and $Z,G/Z$ satisfy bitorsor factorization,
        \item $\Gamma$ is a tree and $G=\SL_1(D)$ where $D$ is a central simple algebra over $F$,
      \end{itemize}
    \item If $\chara(k)=p>0$, assume that $Z(G)$ has finite order not divided by $\chara(k)$ and $G$ and $\Gamma$ are as in the case of $\chara(k)=0$.
    \end{itemize}
  Then, the local-global principle holds for $X$ over $\F$.
\end{corollary} 
\begin{proof}
  Follows from Theorem \ref{curve_local} and Theorem \ref{homspace}.
\end{proof}


\appendix
\section{Nonabelian Hypercohomology}\label{hyper}
Throughout this section, let $F$ be a field and $\Gamma$ its absolute Galois group.
Furthermore, let $H,G$ be a algebraic groups over $F$.
We will review here basics from nonabelian hypercohomology. Nonabelian hypercohomology
was introduced by Breen (\cite{Breen1990}) and Borovoi (\cite{Borovoi1998},\cite{Borovoi1992}).
We refer to these papers for the definition of $\HH^i(F, H \rightarrow G)$ for $i=-1,0,1$.
We will write $\HH^i(G)$ and $\HH^i(G \rightarrow H)$ for $\HH^i(F,G)$ and $\HH^i(F,G \rightarrow H)$ respectively.

\begin{defn}
  A crossed module over $K$ is a morphism $\rho \colon G \rightarrow H$ of algebraic groups over $K$ together with a left action
  $\alpha \colon H \times G \rightarrow G$ such that 
  \begin{itemize}
    \item $\prescript{\rho(g')}{}{g}=g'g (g')^{-1}$  for $g',g \in G$ and
    \item $\rho(\prescript{h}{}{g})=h \rho(g) h^{-1}$ for $g \in G$ and $h \in H$
  \end{itemize}
  holds.
\end{defn}

Given a crossed module $\rho \colon G \rightarrow H$, a $\Gamma$ action on
$\rho$ consists of actions of $\Gamma$ on $G$ and $H$  satisfying
    \begin{align*}
      \rho(\prescript{\sigma}{}{g})=\prescript{\sigma}{}{\rho(g)}
    \end{align*}
and
    \begin{align*}
      \prescript{\sigma}{}{(\prescript{h}{}{g})}=
       \prescript{\prescript{\sigma}{}{h}}{}{(\prescript{\sigma}{}{g})}
    \end{align*}
for $g \in G$, $h \in H$ and $\sigma \in \Gamma$.

\begin{examples}
  We will mostly use the following two examples. 
  \begin{itemize}
    \item Let $G$ be any algebraic group, then $1 \rightarrow G$ is a crossed
      module.
    \item For any algebraic group $G$, $G \xrightarrow{\text{Int}} \Aut(G)$ is
      a crossed module.
  \end{itemize}
\end{examples}

\begin{proposition}[{\cite[Section 3.4.2]{Borovoi1998}}]
\label{hyperles}
  Let  
  \begin{align*}
    1 \rightarrow (G_1 \rightarrow H_1) \xrightarrow{i} (G_2 \rightarrow H_2) \xrightarrow{j}
    (G_3 \rightarrow H_3) \rightarrow 1
  \end{align*}
  be an exact sequence of complexes of $\Gamma$-groups where $i$ is an embedding
  of crossed modules with $\Gamma$ action.
  Then, there is an exact sequence of pointed sets 
  \begin{center}
    \begin{tikzcd}
    1 \arrow{r} &\HH^{-1}(G_1 \rightarrow H_1) \arrow{r} &\HH^{-1}(G_2 \rightarrow H_2) \arrow{r} \arrow[draw=none]{d}[name=Z, shape=coordinate]{} & \HH^{-1}(G_3 \rightarrow H_3)
     \arrow[rounded corners,to path={ -- ([xshift=2ex]\tikztostart.east)|- (Z) [near end]\tikztonodes-| ([xshift=-2ex]\tikztotarget.west)-- (\tikztotarget)}]{dll}\\
    \;& \HH^{0}(G_1 \rightarrow H_1) \arrow{r} &  \HH^{0}(G_2 \rightarrow H_2) \arrow{r} \arrow[draw=none]{d}[name=Z2, shape=coordinate]{} & \HH^{0}(G_3 \rightarrow H_3)
     \arrow[rounded corners,to path={ -- ([xshift=2ex]\tikztostart.east)|- (Z2) [near end]\tikztonodes-| ([xshift=-2ex]\tikztotarget.west)-- (\tikztotarget)}]{dll}\\
   \; & \HH^{1}(G_1 \rightarrow H_1) \arrow{r} & \HH^{1}(G_2 \rightarrow H_2) & \;
  \end{tikzcd}
  \end{center}
\end{proposition}

\begin{example}\label{myles}
  We will mostly use Proposition \ref{hyperles} for the following short exact sequence: 
  \begin{align*}
   1 \rightarrow (1 \rightarrow \Aut(G)) \xrightarrow{i} (G \rightarrow \Aut(G)) \xrightarrow{j}
    (G \rightarrow 1) \rightarrow 1
  \end{align*}
  where all maps occurring are either the identity or trivial. Then, the corresponding long exact sequence simplifies to
  the following sequence (cf. \cite[Section 4.2.3]{Breen1990}).

 \begin{center}
    \begin{tikzcd}
    1 \arrow{r} &\HH^{0}(Z(G)) \arrow{r} &\HH^{0}(G) \arrow{r} \arrow[draw=none]{d}[name=Z, shape=coordinate]{} & \HH^{0}(\Aut(G))
     \arrow[rounded corners,to path={ -- ([xshift=2ex]\tikztostart.east)|- (Z) [near end]\tikztonodes-| ([xshift=-2ex]\tikztotarget.west)-- (\tikztotarget)}]{dll}\\
    \;& \HH^{0}(G \rightarrow \Aut(G)) \arrow{r} &  \HH^{1}(G) \arrow{r} \arrow[draw=none]{d}[name=Z2, shape=coordinate]{} & \HH^{1}(\Aut(G))
     \arrow[rounded corners,to path={ -- ([xshift=2ex]\tikztostart.east)|- (Z2) [near end]\tikztonodes-| ([xshift=-2ex]\tikztotarget.west)-- (\tikztotarget)}]{dll}\\
   \; & \HH^{1}(G \rightarrow \Aut(G))) & \; & \;
  \end{tikzcd}
  \end{center}
\end{example}

We will now describe characterizations of $\HH^i( G \rightarrow \Aut(G))$ for $i=-1,0,1$.
Recall that $\HH^{-1}(H \xrightarrow{\alpha} G)= \ker(\alpha)^{\Gamma}$ and hence
\begin{align*}
  \HH^{-1}(F,G \rightarrow \Aut(G)) = Z(G)(F)
\end{align*}
 where $Z(G)$ is the center of $G$.

As mentioned in section \ref{bitorsor}, $\HH^0(G \rightarrow \Aut(G))$ classifies
$G$-bitorsors. 
\begin{proposition}[{\cite[Theorem 4.5]{Breen1990}}] \label{hyperbit}
 There is a natural isomorphism 
  \begin{align*}
    \HH^0(F,G \rightarrow \Aut(G)) \simeq \left \{ \text{Isomorphism classes of $G$-bitorsors
    over $F$}\right \}.
  \end{align*}
\end{proposition}

Unlike torsors, a $G$-bitorsor may not be trivial (i.e. isomorphic to $G$ as a bitorsor) even when it admits a point. This phenomenon motives the next definition.

\begin{defn}
  Let $\alpha \in \HH^0(G \rightarrow \Aut(G))$. We say that $\alpha$ is \textit{neutral} if
  a bitorsor representing $\alpha$ admits a point over $F$.
\end{defn}

We now turn our attention to $\HH^1(G \rightarrow \Aut(G))$.

\begin{proposition}[{\cite[Theorem 4.5]{Breen1990}}] \label{hypergerbe}
 There is a natural isomorphism 
  \begin{align*}
    \HH^1(F,G \rightarrow \Aut(G)) \simeq \left \{ \text{Equivalence classes of $G$-gerbes
    over $F$}\right \}.
  \end{align*}
\end{proposition}

\begin{defn}
  Let $\alpha \in \HH^1(G \rightarrow \Aut(G))$. We say that $\alpha$ is \textit{neutral} if
  a corresponding gerbe (and thus every corresponding bitorsor) admits a point over $F$.
\end{defn}



\nocite{Harbater2012a}
\nocite{Harbater2011}
\nocite{Harbater2010}
\nocite{Breen}
\nocite{Breen2005}
\nocite{Breen1991}
\nocite{Breen1990}
\nocite{Borovoi1993}
\nocite{Gonzalez-Aviles2012}
\nocite{Borovoi1995}

\bibliographystyle{alpha}
\bibliography{library}

\end{document}